\documentclass{amsart}
\usepackage[utf8]{inputenc}
\usepackage{enumerate}
\usepackage{enumitem}
\usepackage{graphicx}
\usepackage{wrapfig}
\usepackage{subcaption}
\usepackage{listings}
\usepackage{xcolor}
\usepackage[foot]{amsaddr}
\definecolor{mygreen}{RGB}{28,112,30} 
\definecolor{mylilas}{RGB}{170,55,241}
\definecolor{myblue}{RGB}{20,20,230}
\definecolor{myred}{RGB}{200,10,30}
\usepackage{microtype}
\usepackage{multirow}
\usepackage{tikz}
\usepackage{array}
\newcounter{rowcount}
\usepackage{amsmath}
\usepackage{amsfonts}
\usepackage{amssymb}
\usepackage{amscd}
\usepackage{amstext}
\usepackage{amsbsy}
\usepackage{amsopn}
\usepackage{amsthm}
\usepackage{thmtools}
\usepackage{amsxtra}
\usepackage{float}
\usepackage{bbm}
\usepackage{comment}
\usepackage{url}
\usepackage[numbers,square]{natbib}
\usepackage{hyperref}

\newtheoremstyle{mytheoremstyle} 
  {12pt}                    
  {\topsep}                    
  {\slshape}                   
  {}                           
  {\bfseries}                   
  {.}                          
  {.5em}                       
  {}  

\newtheoremstyle{mydefstyle} 
  {12pt}                    
  {12pt}                    
  {}                   
  {}                           
  {\bfseries}                   
  {.}                          
  {.5em}                       
  {}  

  \newtheoremstyle{plainsl}%
  	{}
  	{\topsep}
  	{\slshape} 
  	{}
  	{\normalfont\bfseries}
  	{.}
  	{ }
  	{}

\theoremstyle{plainsl}
\newtheorem{thm}{Theorem}[section]
\newtheorem*{thm*}{Theorem}
\newtheorem{prop}[thm]{Proposition}
\newtheorem{cor}[thm]{Corollary}
\newtheorem{lem}[thm]{Lemma}
\newtheorem*{lem*}{Lemma}
\newtheorem{quest}[thm]{Question}
\newtheorem{prob}[thm]{Problem}

\theoremstyle{remark}
\newtheorem{rem}[thm]{Remark}

\theoremstyle{mydefstyle}

\newtheorem{defn}[thm]{Definition}

\newtheorem{ex}[thm]{Example}
\newtheorem*{note*}{Remark}

\hypersetup{
    plainpages=false,       
    unicode=false,          
    pdftoolbar=true,        
    pdfmenubar=true,        
    pdffitwindow=false,     
    pdfstartview={FitH},    
    pdfnewwindow=true,      
    colorlinks=true,        
    linkcolor=black,         
    citecolor=black,        
    filecolor=blue,      
    urlcolor=blue           
}

\newcommand{\Z}{\mathbb{Z}}

\newcommand{\FF}{\mathbb{F}}

\renewcommand{\phi}{\varphi}

\renewcommand\qed{%
	\ifmmode\eqno\sqr53
	\else\nolinebreak\ \hfill\sqr53\medbreak\fi}
\renewcommand\proof{\noindent\textsl{Proof. }}
\newcommand\sqr[2]{{\vbox{\hrule height.#2pt
    \hbox{\vrule width.#2pt height#1pt \kern#1pt
        \vrule width.#2pt}\hrule height.#2pt}}}

\def\pmat#1{{\begin{pmatrix}#1\end{pmatrix}}}

\newcommand{\qbino}[2]{\left [ #1 \right ]_{#2}}

\DeclareMathOperator\Sym{Sym}
\DeclareMathOperator\Aut{Aut}
\DeclareMathOperator\Cay{Cay}
\DeclareMathOperator\BCay{Cin}
\DeclareMathOperator\HCay{CH}

\DeclareMathOperator\dev{dev}

\DeclareMathOperator\Stab{Stab}
\DeclareMathOperator\dic{Dic}
\DeclareMathOperator\dih{Dih}

\newcommand{\abs}[1]{\left | #1 \right |}
\newcommand{\lp}{\! \left (}
\newcommand{\rp}{\right )}
\newcommand{\lb}{\left [}
\newcommand{\rb}{\right ]}
\newcommand{\lsb}{\left \{ }
\newcommand{\rsb}{\right \} }

\newcommand{\vtxa}{u}
\newcommand{\vtxb}{v}

\newcommand{\bipartC}{\beta}

\newcommand{\bipartB}{\gamma}

\definecolor{mulberry}{rgb}{0.77, 0.29, 0.55}
\definecolor{myred}{RGB}{200,10,30}
\definecolor{mygreen}{RGB}{28,112,30}
\definecolor{UWblack}{RGB}{0,0,0}
\definecolor{UWgold}{RGB}{255,213,79}
\definecolor{anothergreen}{RGB}{120,200,150} 

\definecolor{caBlue}{RGB}{100,143,255}
\definecolor{caPurple}{RGB}{120,94,240}
\definecolor{caRed}{RGB}{220,38,127}
\definecolor{caOrange}{RGB}{254,97,0}
\definecolor{caYellow}{RGB}{255,176,0}

\usepackage{amsaddr}

\title{Cayley incidence graphs}

\author{$^*$Arnbj\"org Soff\'ia \'Arnad\'ottir}
\address{$^*$Universidade Federal de Minas Gerais, Av.\ Pres.\ Antônio Carlos 6627 - Pampulha, Belo Horizonte - MG, 31270-901, Brazil}

\author{$^\dagger$Alexey Gordeev}
\author{$^\dagger$Sabrina Lato}
\author{$^\dagger$Tovohery Randrianarisoa}
\author{$^\dagger$Joannes Vermant}
\address{$^\dagger$Umeå Universitet, Universitetstorget 4, 901 87 Umeå, Sweden}

\email{arnbjorg.soffia@dcc.ufmg.br, aleksei.gordeev@umu.se}
\email{sabrina.lato@umu.se, tovohery.randrianarisoa@umu.se, joannes.vermant@umu.se}

\begin{document}
\renewcommand{\itshape}{\slshape}

\begin{abstract}
Evra, Feigon, Maurischat, and Parzanchevski (2023) introduced a biregular extension of Cayley graphs. In this paper, we reformulate their definition and provide some basic properties. We also show how these Cayley incidence graphs relate to various notions of Cayley hypergraphs. We further establish connections between Cayley incidence graphs and certain geometric and combinatorial structures, including coset geometries, difference sets and cages.

\end{abstract}

\maketitle

\setlength{\parskip}{0cm}
\tableofcontents
\setlength{\parskip}{0.15cm}

\section{Introduction}
Cayley graphs were introduced by Arthur Cayley~\cite{cayley1878} in the late 19th century with the objective of illustrating abstract groups. They have now become a widely studied topic in both group theory and graph theory as they connect these two areas of mathematics in a beautiful way. A Cayley graph is defined by a group and a subset of this group. The group is contained in the automorphism group of the graph and acts regularly on the vertices. In particular, Cayley graphs are vertex-transitive,  which makes them a well-behaved class of graphs to work with.

A key application of Cayley graphs is to explicitly construct graph expanders. Expanders are infinite families of sparse graphs with high connectivity that have a number of applications~\cite{expanderApplications}. They can equivalently be seen as families of graphs where the absolute values of the non-trivial eigenvalues are small. 

For an infinite family of $k$-regular graphs, Alon and Boppana~\cite{alonExpanders} and Serre~\cite{serre} proved that the limit of the sequence of second-largest eigenvalues of graphs in the family is at least $2 \sqrt{k-1}.$ Graphs that meet this bound are called \textit{Ramanujan graphs}, and Lubotzky, Philips, and Sarnak~\cite{lubotzky}, Margulis~\cite{margulis} and Morgenstern~\cite{morgenstern}, among others, used Cayley graphs to explicitly construct $k$-regular Ramanujan graphs for certain values of $k$.

For an infinite family of bipartite biregular graphs, the corresponding bound comes from Feng and Li~\cite{fengLi}. Marcus, Spielman, and Srivastava~\cite{allExpanders} proved for any \( k \) and \( \ell, \) there exists an infinite family of \( \lp \ell, k \rp\)-biregular Ramanujan graphs. Their proof is probablistic, not constructive, and Evra, Feigon, Maurischat, and Parzanchevski~\cite{evra2023ramanujan} introduced what they called Cayley bigraphs to explicitly construct the first bipartite biregular expanders. These Cayley bigraphs are the topic of this paper, but to avoid confusion with Cayley digraphs, we will refer to them as Cayley incidence graphs.

Proving that an infinite sequence of graphs is a family of expanders generally requires deep tools from representation theory and number theory. However, Cayley graphs can also be used in a more sporadic way to construct examples of graphs with certain structural and spectral properties. In Proposition~\ref{bCSpectrum}, we show that Cayley incidence graphs have similar potential by showing that we can compute the spectrum of a Cayley incidence graph at least as easily as we can compute the spectrum of a Cayley graph. Although this is not trivial, computing the spectra of Cayley graphs has been well studied. In 1979, Babai~\cite{babai-spectra} derived an expression for the eigenvalues of a Cayley graph in terms of its irreducible characters. More has been done since then and a survey on the topic was published recently \cite{Liu-survey}.

Any bipartite graph can be interpreted as the incidence graph of a hypergraph, and several different hypergraph generalizations of Cayley graphs have been introduced over the years~\cite{buratti1994cayley, jajcayova2024generalizations, lee2013cayley, shee1990group, teh1976algebraic}. A particular axiom defined by Evra et al~\cite{evra2023ramanujan} imposes additional structure on the incidence graph, and in Corollary~\ref{linearSabidussi}, we show that Cayley incidence graphs are precisely the restriction of the hypergraphs of Shee~\cite{shee1990group}, Lee and Kwon~\cite{lee2013cayley}, or Jajcay and Jajcayov\'a~\cite{jajcayova2024generalizations} to linear uniform regular hypergraphs.

Cayley incidence graphs also relate to structures of interest in design theory and finite geometry.
In Proposition~\ref{prop-diffset-cin} we show that $(n, k, 1)$-difference sets~\cite{beth1999design, jungnickel200618difference} correspond to a certain subclass of Cayley incidence graphs. In particular, it implies that Desarguesian finite projective planes can be represented as Cayley incidence graphs. We show this explicitly by giving constructions of finite projective and affine spaces of any dimension over any finite field as Cayley incidence graphs in Example~\ref{exAGnq} and Example~\ref{exPGnq}. 

Another example comes from the coset geometries introduced by Tits~\cite{tits1963geometries} to characterize groups. In Proposition~\ref{l2} and Example~\ref{cosetGeo}, we illustrate how certain kinds of Cayley incidence graphs relate to a particular class of coset geometries.

If we restrict ourselves to Cayley graphs which are bipartite, a natural question is whether they are related to Cayley incidence graphs. Indeed, as we show in Theorem \ref{cayleyBiCayley}, a bipartite Cayley graph of girth at least six can be represented as a Cayley incidence graph. However, there are some Cayley incidence graphs which cannot be constructed as bipartite Cayley graphs. A generalization of Cayley graphs are bi-Cayley graphs, which admit a semiregular group of automorphisms with two orbits. These graphs have been widely researched and they are significant in studying the symmetry of graphs \cite{Du2015, Du1999151}.
In Theorem \ref{BicCayG00S} we show that a bipartite bi-Cayley graph with girth at least six is isomorphic to a Cayley incidence graph.

Cayley incidence graphs are a new class of graphs with close connections to generalizations of Cayley graphs studied in various contexts. We conclude the paper by highlighting some open problems of particular interest, though much of the theory and applications of these graphs remains open for further study.

\section{Definition}

Throughout this paper, we will generally assume that \( G \) is a multiplicative group with identity \( e \). However, when dealing with specific groups it will often be convenient to use the notation for that group as, for instance, using additive notation and identity 0 when dealing with cyclic groups.

For a set \( C \) of elements of \( G, \) we will let 
\[ gC := \lsb gs: s \in C \rsb. \]

Let \( S \) be an inverse-closed subset of \( G \setminus \lsb e \rsb. \) The \textit{Cayley graph}, \( \Cay \lp G, S \rp \), is the graph on vertex set \( G \) with vertices \( g, h \) adjacent if and only if \( h^{-1} g \in S. \) The set $S$ is called the \emph{connection set} of $\Cay(G,S)$.

Let $G$ be a group, and let $\pi = \lsb C_1, \ldots, C_{\ell} \rsb$ be a collection of subsets of $G$ that each contains the identity.
Let 
\[ \bipartB := G \text{ and } \bipartC := \lsb gC : g \in G, C \in \pi \rsb. \]
We wish to define an incidence relation between \( \bipartB \) and \( \bipartC \). If \( g \in \bipartB \), we can define an incidence relation \( g \sim_1 gC_i \) for \( 1 \leq i \leq \ell. \) If \( gC \in \bipartC, \) we can define an incidence relation \( gC \sim_2 gs \) for all \( s \in C. \) Both these incidence relations are natural, but they are not generally equivalent.

\begin{ex}Let \( G = \mathbb{Z}_5 \) and let \( \pi = \lsb \lsb 0, 1, 2 \rsb, \lsb 0, 3, 4 \rsb \rsb. \) Using \( \sim_1, \) we see that \( 1 \) is incident to \( \lsb 1, 2, 3 \rsb \) and \( \lsb 0, 1, 4 \rsb. \) However, using \( \sim_2, \) we have that 1 is incident to \( \lsb 0, 1, 2 \rsb, \lsb 1, 2, 3 \rsb \) and \( \lsb 0, 1, 4 \rsb. \)
\end{ex}

Fix \( D \in \bipartC. \) It is incident to the same set of elements from \( G \) due to $\sim_1$ and $\sim_2$ precisely when, for every $h \in D,$ there exists $1 \leq j \leq \ell$ such that $h^{-1}D = C_j$. We can reframe this as a condition on \( \pi \).

\begin{defn}[$T$-axiom]\label{abc}
    Let \( G \) be a group and let \( \pi \) be a collection of subsets of $G$ that each contain the identity. We say it satisfies the \textit{$T$-axiom} if for any $C \in \pi$ and $s \in C,$ we have $s^{-1} C \in \pi$.
\end{defn}

The $T$ in \nameref{abc} stands for ``translate'' as it imposes a condition on translates of $C \in \pi$.
We will refer to the sets $C \in \pi$ as \textit{cells}. We can now use our singular incidence relation to define a bipartite incidence graph between group elements and translates of the cells.

\begin{defn}[Cayley incidence graph]\label{bcayDef}
Let $G$ be a group, let $\pi = \lsb C_1, \ldots, C_{\ell} \rsb$ be a collection of subsets of $G$ of size $k$ that intersect pairwise precisely in the identity and suppose $\pi$ satisfies the \nameref{abc}. Let
    \[ \bipartB = G \text{ and } \bipartC = \lsb gC : g \in G, C \in \pi \rsb. \]
    We define the \textit{Cayley incidence graph} \( \BCay \lp G, \pi \rp \) as the bipartite graph on vertex set \( \bipartB \cup \bipartC \) with an edge from \( g \) to \( gC \) for all \( g \in G \) and \( C \in \pi\).
\end{defn} 

Note that \( \BCay \lp G,  \pi \rp \) is biregular. An \( \lp \ell, k \rp \)-biregular Cayley incidence graph has $\ell$ cells of size $k$. 

\begin{ex}\label{exCoarse}If \( \ell =1, \) then \( \pi \) has a single cell \( C \). This satisfies the \nameref{abc} precisely when \( C \) is a subgroup. In this case \( \BCay \lp G, \lsb C \rsb \rp \) is a disjoint union of \( K_{1,k} \).
\end{ex}

\begin{ex}\label{exSubgroups} More generally, whenever $\pi$ is a collection of $\ell$ subgroups of $G$ of order $k$ that intersect pairwise in the identity, the \nameref{abc} is trivially satisfied, so $\BCay(G, \pi)$ is well-defined.
\end{ex}

For any $\pi$ coming from a Cayley incidence graph, let
\[ S \lp \pi \rp = \bigcup_{i=1}^{\ell} C_i \setminus \lsb e \rsb. \]
It is easy to see that $S(\pi)$ is an inverse-closed subset of $G\setminus \{e\}$, so we may define the Cayley graph $\Cay(G,S(\pi))$. This is the \emph{underlying Cayley graph} of $\BCay(G,\pi)$.

\begin{ex}\label{exDiscrete}If \( k =2,\) then \( \pi = \lsb \lsb e, g_1 \rsb, \lsb e, g_2 \rsb, \ldots, \lsb e, g_{\ell} \rsb \rsb. \) Then \( \pi \) satisfies the \nameref{abc} precisely when \( S \lp \pi \rp \) is inverse-closed. In this case, \( \BCay \lp G, \pi \rp \) is the subdivision graph of the underlying Cayley graph \( \Cay \lp G, S \lp \pi \rp\rp. \)
\end{ex}

Note that any $gC \in \bipartC$ corresponds to a clique in the underlying Cayley graph: if $gh_1, gh_2 \in gC$, then by the \nameref{abc} $h_2^{-1}C \in \pi$, so $h_2^{-1}g^{-1}gh_1 = h_2^{-1}h_1 \in h_2^{-1}C \subseteq S(\pi)$.
In particular, $\pi$ corresponds to a set of $\ell$ cliques of size $k$ in $\Cay(G, S(\pi))$ pairwise intersecting in the identity.

The cases where \( \ell =1 \) or \( k =2 \) are trivial, and we will generally be interested in the non-trivial examples with  \( \ell \geq 2 \) and \( k \geq 3.\)

\begin{ex}[Fano Plane]\label{Z7}Let \( G = \mathbb{Z}_7 \) and \( \pi = \lsb \lsb 0, 1, 3 \rsb, \lsb 0, 2, 6 \rsb, \lsb 0, 4, 5 \rsb \rsb. \) We can verify that this satisfies the \nameref{abc}. The resulting graph, shown in Figure~\ref{heawood}, is the smallest non-trivial Cayley incidence graph. Sometimes called the \textit{Heawood graph}, it can also be viewed as the incidence graph of the Fano plane. \end{ex}

\begin{figure}[h]
\begin{minipage}{2in}
  \begin{center}
    \begin{tikzpicture}[node distance = 2 cm]
      \tikzset{VertexStyle/.style = {shape= circle, ball color=UWblack, inner sep= 2pt, outer sep= 0pt, minimum size = 24 pt, text= caYellow}}
      \tikzset{EdgeStyle/.style   = {double= UWblack, double distance = 1pt}}
    
      \draw (90:2.5 cm) node[VertexStyle, font=\small] (0) {0} ;
      \draw (141.4:2.5 cm) node[VertexStyle, font=\small] (12) {4} ;
      \draw (192.8:2.5 cm) node[VertexStyle, font=\small] (10) {3} ;
      \draw (244.2:2.5 cm) node[VertexStyle, font=\small] (8) {5} ;
      \draw (295.6:2.5 cm) node[VertexStyle, font=\small] (6) {6} ;
      \draw (347:2.5 cm) node[VertexStyle, font=\small] (4) {2} ;
      \draw (39.4:2.5 cm) node[VertexStyle, font=\small] (2) {1} ;

      \tikzset{VertexStyle/.style = {shape= circle, ball color= caYellow!60!white, inner sep= 2pt, outer sep= 0pt, minimum size = 24 pt, text=UWblack}}

      \draw (115.7:2.5 cm) node[VertexStyle, font=\scriptsize] (13) {$\lsb 4, 5, 0 \rsb$} ;
      \draw (167.1:2.5 cm) node[VertexStyle, font=\scriptsize] (11) {$\lsb 3, 4, 6 \rsb$} ;
      \draw (218.5:2.5 cm) node[VertexStyle, font=\scriptsize] (9) {$\lsb 2, 3, 5 \rsb$} ;
      \draw (269.9:2.5 cm) node[VertexStyle, font=\scriptsize] (7) {$\lsb 5, 6, 1 \rsb$} ;
      \draw (321.3:2.5 cm) node[VertexStyle, font=\scriptsize] (5) {$\lsb 6, 0, 2 \rsb$} ;
      \draw (12.7:2.5 cm) node[VertexStyle, font=\scriptsize] (3) {$\lsb 1, 2, 4 \rsb$} ;
      \draw (65.1:2.5 cm) node[VertexStyle, font=\scriptsize] (1) {$\lsb 0, 1, 3 \rsb$} ;

      \draw [EdgeStyle] (0)--(1) ;
      \draw [EdgeStyle] (1)--(2) ;
      \draw [EdgeStyle] (2)--(3) ;
      \draw [EdgeStyle] (3)--(4) ;
      \draw [EdgeStyle] (4)--(5) ;
      \draw [EdgeStyle] (5)--(6) ;
      \draw [EdgeStyle] (6)--(7) ;
      \draw [EdgeStyle] (7)--(8) ;
      \draw [EdgeStyle] (8)--(9) ;
      \draw [EdgeStyle] (9)--(10) ;
      \draw [EdgeStyle] (10)--(11) ;
      \draw [EdgeStyle] (11)--(12) ;
      \draw [EdgeStyle] (12)--(13) ;
      \draw [EdgeStyle] (13)--(0) ;

      \draw [EdgeStyle] (0)--(5) ;
      \draw [EdgeStyle] (1)--(10) ;
      \draw [EdgeStyle] (2)--(7) ;
      \draw [EdgeStyle] (3)--(12) ;
      \draw [EdgeStyle] (4)--(9) ;
      \draw [EdgeStyle] (6)--(11) ;
      \draw [EdgeStyle] (8)--(13) ;
         
    \end{tikzpicture}
  \end{center}
\end{minipage}
\hspace{.5in}
\begin{minipage}{2in}
  \begin{center}
    \begin{tikzpicture}[node distance = 2 cm]
      \tikzset{VertexStyle/.style = {shape= circle, ball color= caYellow!60!white, inner sep= 2pt, outer sep= 0pt, minimum size = 24 pt, text= UWblack}}
      \tikzset{EdgeStyle/.style   = {double= UWblack, double distance = 1pt}}
    
      \draw (90:3 cm) node[VertexStyle, font=\small] (0) {0} ;
      \draw (210:3 cm) node[VertexStyle, font=\small] (1) {2} ;
      \draw (330:3 cm) node[VertexStyle, font=\small] (2) {5} ;
      \draw (270:1.6 cm) node[VertexStyle, font=\small] (3) {3} ;
      \draw (30:1.6 cm) node[VertexStyle, font=\small] (4) {4} ;
      \draw (150:1.6 cm) node[VertexStyle, font=\small] (5) {6} ;      
      \draw (0:0 cm) node[VertexStyle, font=\small] (6) {1} ;

      \draw [EdgeStyle] (0)--(5) ;
      \draw [EdgeStyle] (5)--(1) ;
      \draw [EdgeStyle] (1)--(3) ;
      \draw [EdgeStyle] (3)--(2) ;
      \draw [EdgeStyle] (2)--(4) ;      
      \draw [EdgeStyle] (4)--(0) ;

      \draw [EdgeStyle] (1)--(6) ;      
      \draw [EdgeStyle] (6)--(4) ;            
      \draw [EdgeStyle] (5)--(6) ;      
      \draw [EdgeStyle] (6)--(2) ;            
      \draw [EdgeStyle] (0)--(6) ;            
      \draw [EdgeStyle] (6)--(3) ;           
      \draw [EdgeStyle] (5) edge [ bend left=40] (4) ;
      \draw [EdgeStyle] (4) edge [ bend left=40] (3) ;      
      \draw [EdgeStyle] (3) edge [ bend left=40] (5) ;

    \end{tikzpicture}
  \end{center}
\end{minipage}
\caption{$\BCay \lp \mathbb{Z}_7, \lsb 0, 1, 3 \rsb, \lsb 0, 2, 6 \rsb, \lsb 0, 4, 5 \rsb \rp$ (Heawood graph) and the Fano plane}\label{heawood} 
\end{figure}

Analogously to Cayley graphs, a Cayley incidence graph $\BCay \lp G, \pi \rp$ is connected if and only if $S \lp \pi \rp$ generates $G$. If $S \lp \pi \rp$ does not generate $G$, then we can without loss of generality let $G'$ be the group generated by $S \lp \pi \rp$ and consider the connected graph $\BCay \lp G', \pi \rp$, which is isomorphic to each connected component of $\BCay \lp G, \pi \rp.$ Thus, when we refer to non-trivial Cayley incidence graphs, we assume that they are connected.

Let $G$ be a group, and let $\pi = \lsb C_1, \ldots, C_{\ell} \rsb$ be a collection of subsets of $G$ such that for $1 \leq i \leq \ell$ the cell $C_i$ has size $k$, contains the identity, and any two distinct cells intersect in precisely the identity. Suppose further that $\pi$ satisfies the \nameref{abc} and $k \leq \ell$. We can define a partition $\hat{\pi}$ of $S \subseteq G$ by removing the identity from every cell in in $\pi$. Such a partition satisfies:
\begin{enumerate}[label=(\arabic*)]
\item\label{bCaxiom1}The set \( S \) is inverse-closed and does not contain the identity and
\item\label{bCaxiom3}Whenever $s$ and $t$ are in the same cell $C \in \hat{\pi},$ the product $s^{-1}t$ is in the cell containing $s^{-1}.$
\end{enumerate}

Conversely, given a group $G$ and a partition $\hat{\pi}$ into $\ell$ cells of size $k-1$ satisfying~\ref{bCaxiom1} and~\ref{bCaxiom3}, we can add the identity to every cell $C \in \hat{\pi}$ to get collection \( \pi \) satisfying the \nameref{abc}. 
Evra, Feigon, Maurischat, and Parzanchevski~\cite{evra2023ramanujan} defined the bi-Cayley axioms as~\ref{bCaxiom1} and~\ref{bCaxiom3}, but we find it more convenient to work with the collection \( \pi \) where every cell contains the identity and use the equivalent formulation of Definition~\ref{abc}. We call this the \nameref{abc} to avoid confusion with the bi-Cayley graphs discussed in Section~\ref{secBiBi}.

Evra, Feigon, Maurischat, and Parzanchevski defined an equivalence relation on \( G \times \hat{\pi} \) and used this to define a bipartite graph on \( \bipartB \cup \bipartC \) with \( \bipartB = G \) and \( \bipartC \) the equivalence classes of their relation. They referred to the resulting graphs as \textit{Cayley bigraphs}, and they are precisely the Cayley incidence graphs with the additional assumption that \( k \leq \ell. \)

Another perspective is to view a Cayley incidence graph $\BCay(G, \pi)$ as the incidence graph of a hypergraph. The corresponding hypergraph $\HCay \lp G, \pi \rp$ is a $k$-uniform, $\ell$-regular hypergraph on vertex set \( \bipartB = G \) with hyperedge set
\[ \bipartC = \lsb gC : g \in G, C \in \pi \rsb. \]

Any graph can be viewed as a 2-uniform hypergraph, and, as we have seen in Example~\ref{exDiscrete}, a Cayley graph $\Cay \lp G, S \rp$ is precisely the 2-uniform hypergraph $\HCay \lp G, \lsb \lsb e, s \rsb: s \in S \rsb \rp.$ Thus viewing Cayley incidence graphs as hypergraphs gives a hypergraph generalization of Cayley graphs. 

This is far from the first hypergraph generalization of Cayley graphs~\cite{buratti1994cayley,jajcayova2024generalizations,lee2013cayley,shee1990group}, and we will explore the varying definitions and their relationship to Cayley incidence graphs in Section~\ref{secHyper}. First, it will be convenient to set up more properties of Cayley incidence graphs. 

\section{Basic Properties}

Let \( X = \lp \bipartB \cup \bipartC, E \rp \) be any bipartite graph.
The \textit{biadjacency matrix} \( N \) of $X$ is the \( \bipartB \times \bipartC \) matrix with \( \lp \vtxa, \vtxb \rp\) entry equal to 1 if \( \vtxb \) is incident to \( \vtxa \) and 0 otherwise. Then the adjacency matrix has the form
\[ \pmat{\mathbf{0} & N \\ N^T & \mathbf{0}}. \]

Consider the distance-two graph \( X_2\) of \( X. \) This is the graph on the vertex set of \( X \) with two vertices adjacent in \( X_2\) if and only if they are at distance two in \( X. \) Since \( X \) is bipartite, the distance-two graph can be written as a disjoint union of a graph on vertex set \( \bipartB \) and a graph on vertex set \( \bipartC. \) These graphs are the \textit{halved graphs}, denoted \( H_{\bipartB} \) and \( H_{\bipartC} \). 

In general, the adjacency matrix of \( H_{\bipartB} \) is given by the matrix with 1 on every nonzero off-diagonal entry of \( NN^T \) and  0 everywhere else. The next proposition tells us that for a Cayley incidence graph, the halved graph \( H_{\bipartB} \) is the underlying Cayley graph, with a particularly strong relationship between the adjacency matrix of $H_{\bipartB}$ and $NN^T.$

\begin{prop}\label{cayleyHalved}Let \( \BCay \lp G, \pi \rp \) be an \( \lp \ell, k \rp \)-biregular Cayley incidence graph with biadjacency matrix $N$ and let $A$ be the adjacency matrix of the underlying Cayley graph \( \Cay \lp G, S \lp \pi \rp \rp \).
Then
\[ NN^T = A + \ell I. \]
\end{prop}

\proof Observe that \( NN^T \) counts the number of walks of length two between elements of \( G. \) Thus, if \( g, h \in G, \) we have that
\[ \lp N N^T \rp_{g,h} = \abs{\lsb gC_i : h \in gC_i \rsb} = \abs{\lsb C_i : g^{-1}h \in C_i \rsb} = 
\begin{cases}
\ell &\text{if } h=g,\\ 
1 & \text{if } g^{-1}h \in S(\pi),\\ 
0 & \text{if }g^{-1}h \notin S(\pi). 
\end{cases}\]
Recalling the definition of Cayley graphs and translating back into matrices gives the desired result.\qed

We defined the underlying Cayley graph in terms of the group and the connection set. Proposition~\ref{cayleyHalved} gives another way of interpreting the relationship between the Cayley incidence graph and its underlying Cayley graph, which is particularly useful in computing the spectrum. 

We will need the following fact from linear algebra, which can be found, for instance, in Section~10.3 of~\cite{blue}. 

\begin{lem}\label{BBT}Let $N$ be a matrix. Then the nonzero eigenvalues of $N^TN$ and $NN^T$ are the same, with the same multiplicity.\qed\end{lem}

\begin{prop}\label{bCSpectrum}The spectrum of a Cayley incidence graph is determined by the spectrum of the underlying Cayley graph, the order of the group, and the valencies $k$ and $\ell$.
\end{prop}

\proof Let $A$ be the adjacency matrix of a Cayley incidence graph $\BCay \lp G, \pi \rp$ and let $N$ be the biadjacency matrix. We have
\[ A^2 = \pmat{NN^T & \mathbf{0} \\ \mathbf{0} & N^TN}. \]

Let the spectrum of the underlying Cayley graph \( \Cay \lp G, S \lp \pi \rp \rp \) be
\[ \lsb \theta_1^{(m_1)}, \ldots, \theta_{t-1}^{(m_{t-1})}, \theta_t^{(m_t)} \rsb \]
where \( \theta_r\) is an eigenvalue and \( m_r \) is the multiplicity of \( \theta_r\) for all \( 1 \leq r \leq t. \) By Proposition~\ref{cayleyHalved}, we see that
\[ \lsb \lp \theta_1+ \ell \rp^{(m_1)}, \lp \theta_2+\ell \rp^{(m_2)}, \ldots, \lp \theta_t+\ell \rp^{(m_t)} \rsb \]
are the eigenvalues of \( NN^T.\)

By Lemma~\ref{BBT}, we know that $NN^T$ and $N^T N$ share nonzero eigenvalues with multiplicity. In particular, if $\theta_r+\ell$ is a nonzero eigenvalue of $NN^T$ with multiplicity $m_r,$ then it is an eigenvalue of $A^2$ with multiplicity $2 m_r.$ All that remains is to determine the multiplicity of zero as an eigenvalue for $A^2.$

If $n$ is the order of $G$, then $\bipartB$ has size $n$ and $\bipartC$ has size $b = n\ell / k$. We suppose for notational convenience that $\theta_t = -\ell,$ possibly with multiplicity $m_t = 0.$ The multiplicity of 0 in $N^TN$ is $m_t+b-n$, so the multiplicity of 0 as an eigenvalue in $A^2$ is $2m_r+b-n$.

Finally, the spectrum of a bipartite graph is symmetric over the real axis, so we get that the spectrum of $\BCay \lp G, \pi \rp$ is
\[\lsb \pm\sqrt{\theta_1+\ell  }^{(m_1)}, \ldots, \pm\sqrt{\theta_{t-1}+\ell}^{(m_{t-1})}, 0^{(2m_t+b-n)} \rsb.\qed\]

The \textit{girth} of a graph is the length of a shortest cycle in it. Proposition~\ref{cayleyHalved} gives us information about the girth of a Cayley incidence graph.

\begin{lem}\label{girth6Plus}A Cayley incidence graph has girth at least six.
\end{lem}

\proof Let $\vtxa, \vtxb \in G$, and recall that $\lp NN^T \rp_{\vtxa, \vtxb}$ counts the number of walks of length two from \( \vtxa \) to \( \vtxb. \) Then for $\vtxa \neq \vtxb,$ we see from Proposition~\ref{cayleyHalved} that the number of paths from \( \vtxa \) to \( \vtxb \) is at most one. In particular, there cannot be a cycle of length four that goes through \( \vtxa, \) and since this is independent of our choice of \( \vtxa, \) there cannot be a cycle of length four in the graph.\qed

Lemma \ref{girth6Plus} implies that every Cayley incidence graph is the incidence graph of a combinatorial configuration, as defined in \cite{grunbaum2009}.

Given two Cayley graphs, $\Cay(G_1,S_1)$ and $\Cay(G_2,S_2)$, their Cartesian product is a Cayley graph for the group $G_1\times G_2$ with connection set $S_1\times \{e\}\cup \{e\}\times S_2$. Thus, there is a natural way of constructing a Cayley graph given two smaller ones. Similarly, we can construct Cayley incidence graphs from smaller Cayley incidence graphs.

\begin{lem}\label{productCIN}Suppose that \( \BCay \lp G_1, \pi_1 \rp \) is an \( \lp \ell_1, k \rp \)-biregular Cayley incidence graph and \( \BCay \lp G_2, \pi_2 \rp \) is an \( \lp \ell_2, k \rp \)-biregular Cayley incidence graph.
We define
    \[ \pi := \{\{e\}\times C\colon C\in \pi_1\}\cup \{C\times \{e\}\colon C\in \pi_2\}. \]
    Then $(G_1 \times G_2, \pi )$ satisfies the \nameref{abc} and $\BCay \lp G_1 \times G_2, \pi \rp = \lp \bipartB \cup \bipartC, E \rp$ has edges of one of the following forms:
    \begin{itemize}
        \item either $\{(g_1,h),g_2C_i\times \{h\}\}$ where $\{g_1,g_2C_i\}$ is an edge of $\BCay(G_1, \pi_1)$, 
        \item or $\{(h,g_1),\{h\}\times g_2C_i]\}$ where $\{g_1,g_2C_i\}$ is an edge of $\BCay(G_2, \pi_2)$.
    \end{itemize}
\end{lem}

\proof It is clear by construction that $\pi$ satisfies the \nameref{abc}.  Now, suppose that in $\BCay(G_1 \times G_2, \pi)$ we have an edge between $(g_1,h_1)$ and $g_2C_i \times \{h_2\}$. This statement is equivalent to $h_1=h_2$ and $g_1 \in g_2C_i,$ which is equivalent to the fact that in the graph $\BCay(G,\pi_1)$, we have an edge between $g_1$ and $g_2C_i$. An analogous result holds for cells of the form $\{e\}\times C_i$.\qed

Note that the underlying Cayley graph of this Cayley incidence graph is the Cartesian product of the underlying Cayley graphs of $\BCay(G_1,\pi_1)$ and $\BCay(G_2,\pi_2)$.

In fact, Lemma~\ref{productCIN} can let us construct non-trivial Cayley incidence graphs from trivial Cayley incidence graphs.

\begin{ex}\label{Z3d}
     Let $G_1$ and $G_2$ be two groups of the same order, possibly isomorphic, and consider the operation in Lemma~\ref{productCIN} applied to $\BCay \lp G_1, \{G_1\} \rp$ and $\BCay \lp G_2, \{G_2\} \rp$. This gives us a Cayley incidence graph on $G_1 \times G_2$ with 
     \[ \pi = \lsb \lp e, g \rp : g \in G_1 \rsb \cup \lsb \lp g, e \rp : g \in G_2 \rsb. \]
     Note that this is an instance of Example~\ref{exSubgroups}.
\end{ex} 

Another way of constructing Cayley incidence graphs comes from taking intersections of two partitions. 

\begin{lem}\label{lem:intersection}
    Let $G$ be a group and let $\pi_1=\{C_1,\dots,C_\ell\}$ and $\pi_2=\{C'_1,\dots,C'_{\ell'}\}$ be collections of subsets of $G$ satisfying the \nameref{abc}. Then the collection 
    \[\pi_1\cap\pi_2:=\{C\cap C': C\in \pi_1, C'\in \pi_2\}\]
    satisfies the \nameref{abc}.
\end{lem}
\proof 
Clearly, each set in $\pi:=\pi_1\cap\pi_2$ contains the identity. Let 
$D\in \pi$ and $s\in D$. Then $D=C\cap C'$ for some $C\in \pi_1$ and $C'\in \pi_2$ and we have $s^{-1}C\in \pi_1$ and $s^{-1}C'\in \pi_2$. Therefore, $s^{-1}D = s^{-1}(C\cap C') = s^{-1}C\cap s^{-1}C'\in \pi$
as desired.\qed

\begin{cor}
    Let $\BCay(G,\pi)$ and $\BCay(G,\pi')$ be two Cayley incidence graphs of $G$. For $r\geq 1,$ define $\pi_r$ by 
    \[\pi_r = \{C_i\cap C_j': C_i\in \pi, C_j'\in \pi', \lvert C_i\cap C_j' \rvert = r\}.\]
    If $\pi_r$ is non-empty, then $\BCay(G,\pi_r)$ is a well-defined Cayley incidence graph of $G$.
\end{cor}
\proof 
Clearly, $\pi_r$ is a collection of subsets of $G$ of a fixed size $r$ and they pairwise intersect precisely in the identity. By Lemma \ref{lem:intersection} and the fact that $\lvert xC\rvert = \lvert C\rvert$, for all $C\subseteq G$ and $x\in G$, we see that each non-empty $\pi_r$ satisfies the \nameref{abc}. Consequently, $\BCay(G,\pi_r)$ is a well-defined Cayley incidence graph.\qed

\begin{ex}
    Let $G=\Z_{15}$ and define 
    \begin{align*}
        \pi &= \left\{\{0,1,4,6\}, \{0,2,11,12\}, \{0,3,5,14\}, \{0,9,10,13\}\right\},\\
        \pi' &= \left\{\{0,1,9,13\}, \{0,2,3,11\}, \{0,4,6,7\}, \{0,8,12,14\}\right\}.\\
    \end{align*}
    We see that $\pi_3 = \left\{\{0,2,11\}, \{0,4,6\},\{0,9,13\}\right\}$ and $\BCay(G,\pi_3)$ is a $3$-regular Cayley incidence graph of $G$. 
\end{ex}

\section{Cayley Hypergraphs}\label{secHyper}

An important result in the theory of Cayley graphs is the following characterization of Sabidussi~\cite{sabidussi1964vertex}.

\begin{thm}[Sabidussi~\cite{sabidussi1964vertex}]\label{sabidussi}A graph \( X \) is a Cayley graph if and only if it admits a group of automorphisms \( G \) acting regularly on the vertex set \( V \lp X \rp. \)
\end{thm}

This characterization motivates many of the definitions of Cayley hypergraphs that have been introduced over the years.

Teh and Shee~\cite{teh1976algebraic} introduced the following hypergraph extension of Cayley graphs, called group hypergraphs. The definition is equivalent to the more recent reintroductions by Shee~\cite{shee1990group} and Lee and Kwon~\cite{lee2013cayley}.

\begin{defn}[Group Hypergraph]\label{groupHypergraph}
    Let \( G \) be a group and let \( \mathcal{C} \) be a collection of non-empty subsets of \( G \). The \textit{group hypergraph} \( \HCay \lp G, \mathcal{C} \rp \) is the hypergraph with vertex set \( G \) and hyperedge set 
\[ \lsb gC : g \in G, C \in \mathcal{C} \rsb. \]
\end{defn}

Note that our Cayley incidence graphs \( \BCay \lp G, \pi \rp \) are the incidence graphs of group hypergraphs \( CH \lp G, \pi \rp \) with substantially more conditions placed on \( \pi. \) 

Analogously to Theorem~\ref{sabidussi}, Lee and Kwon~\cite{lee2013cayley} proved that a hypergraph is a group hypergraph if and only it admits a group of automorphisms acting regularly on the vertex set. A similar result in the context of vertex transitive incidence structures was proven previously by Jajcay~\cite{jajcay2002representing}.

This led Jajcay and Jajcayov\'a~\cite{jajcayova2024generalizations} to define the following version of Cayley hypergraphs.

\begin{defn}[Uniform Cayley hypergraph]\label{uniformCayley}
A \textit{uniform Cayley hypergraph} on vertex set $V$ is an $r$-uniform hypergraph for some $1 \leq r \leq \abs{V}$ that admits a group of automorphisms acting regularly on its vertices.
\end{defn}

A more constructive version of Definition~\ref{uniformCayley} comes from the restriction of Definition~\ref{groupHypergraph} to $r$-uniform hypergraphs.
We can go one step further and ask for a version of Theorem~\ref{hyperSabidussi} for uniform and regular hypergraphs.

\begin{thm}\label{hyperSabidussi}
    A hypergraph $X$ is $k$-uniform and $\ell$-regular with a group $G$ acting regularly on the vertex set if and only if there exists a collection $\mathcal{C} = \lsb C_1, \ldots, C_{\ell} \rsb$ of subsets of $G$ of size $k$ that each contains the identity and satisfies the \nameref{abc}, in which case $X$ is isomorphic to $CH \lp G, \mathcal{C} \rp.$    
\end{thm}
\proof We have already seen that all hypergraphs defined using a partition $\pi$ satisfying the \nameref{abc} are $k$-uniform and $\ell$-regular. We may thus suppose $X$ is a hypergraph admitting a group action which is regular on the vertex set. 

By the variants of Sabidussi's theorem proved in \cite{jajcay2002representing, lee2013cayley}, we know that $X$ is isomorphic to a group hypergraph \( CH \lp G, \mathcal{C} \rp \) for some collection $\mathcal{C}$. We may assume without loss of generality that $\mathcal{C}$ is the set of hyperedges that contain the identity.

    Since \( X \) is $\ell$-regular, we know that \( \abs{\mathcal{C}} = \ell,\)
and since \( X \) is $k$-uniform, every set $C \in \mathcal{C}$ must have size $k$. Finally, if $g \in C$ for some $C \in \mathcal{C},$ then $e \in g^{-1} C,$ so $g^{-1}C \in \mathcal{C}$ and thus $\mathcal{C}$ satisfies the \nameref{abc}.\qed

A hypergraph is \textit{linear} if any two hyperedges intersect in at most one vertex or, equivalently, if its incidence graph has no four-cycles. From Lemma~\ref{girth6Plus}, we see that the hypergraphs arising from Cayley incidence graphs are linear. In fact, Cayley incidence graphs are precisely the linear hypergraphs of Theorem~\ref{hyperSabidussi}.

\begin{cor}\label{linearSabidussi}
    A linear hypergraph is uniform and regular with a group acting regularly on the vertex set if and only if its incidence graph is a Cayley incidence graph.
\end{cor}

Theorem~\ref{hyperSabidussi} begs the question of why we would restrict ourselves to linear hypergraphs when, by placing less conditions on \( \pi, \) we get a larger class of uniform, regular hypergraphs satisfying a hypergraph version of Theorem~\ref{sabidussi}.

Let \( \mathcal{C} = \lsb C_1, \ldots, C_{\ell} \rsb \) be a collection of subsets of \( G \) such for all \( 1 \leq i \leq \ell, \) the cell \( C_i \) has size $k$ and contains the identity. If \( \mathcal{C} \) satisfies the \nameref{abc}, we can define \( X \) to be the bipartite incidence graph of \( CH \lp G, \mathcal{C} \rp\), that is, \( X \) has bipartition \( \bipartB = G\) and
\[ \bipartC = \lsb gC : g \in G, C \in \mathcal{C} \rsb. \]

Note that 
\[ S \lp \mathcal{C} \rp = \bigcup_{C \in \mathcal{C}} C \setminus \lsb e \rsb \]
is still inverse-closed, so we can still speak about the underlying Cayley graph. However, the relationship between \( X \) and this underlying Cayley graph is weaker than the relationship for the Cayley incidence graphs. In the proof of Proposition~\ref{cayleyHalved}, the key fact was that that \( g^{-1}h \neq e \) appears in a single cell of \( \pi \) if \( g^{-1} h \in S \) and in no cells of \( \pi \) otherwise. For \( \mathcal{C} \), there can be vertices in \( \bipartB \) at distance two with a differing number of walks of length two between them. This has particularly devastating consequences for finding the spectrum, since we cannot determine it from the underlying Cayley graph.

We can remedy this problem by further requiring that every non-identity element that appears in $\mathcal{C}$ appears in some constant $\lambda$ number of cells. However, this nicer linear algebraic structure takes away from the group theoretic structure.

\begin{figure}
\begin{center}
    \begin{tikzpicture}[node distance = 2 cm]
      \tikzset{VertexStyle/.style = {shape= circle, ball color= UWblack, inner sep= 2pt, outer sep= 0pt, minimum size = 24 pt, text= caYellow}}
      \tikzset{EdgeStyle/.style   = {double= UWblack, double distance = 1pt}}
    
      \draw (90:3 cm) node[VertexStyle, font=\small] (0) {0} ;
      \draw (141.4:3 cm) node[VertexStyle, font=\small] (12) {4} ;
      \draw (192.8:3 cm) node[VertexStyle, font=\small] (10) {3} ;
      \draw (244.2:3 cm) node[VertexStyle, font=\small] (8) {5} ;
      \draw (295.6:3 cm) node[VertexStyle, font=\small] (6) {6} ;
      \draw (347:3 cm) node[VertexStyle, font=\small] (4) {2} ;
      \draw (39.4:3 cm) node[VertexStyle, font=\small] (2) {1} ;

      \tikzset{VertexStyle/.style = {shape= circle, ball color= caYellow!60!white, inner sep= 2pt, outer sep= 0pt, minimum size = 24 pt, text=UWblack}}

      \draw (115.7:3 cm) node[VertexStyle, font=\scriptsize] (13) {$\lsb 5, 6, 0 \rsb$} ;
      \draw (167.1:3 cm) node[VertexStyle, font=\scriptsize] (11) {$\lsb 3, 4, 5 \rsb$} ;
      \draw (218.5:3 cm) node[VertexStyle, font=\scriptsize] (9) {$\lsb 2, 3, 4 \rsb$} ;
      \draw (269.9:3 cm) node[VertexStyle, font=\scriptsize] (7) {$\lsb 4, 5, 6 \rsb$} ;
      \draw (321.3:3 cm) node[VertexStyle, font=\scriptsize] (5) {$\lsb 6, 0, 1 \rsb$} ;
      \draw (12.7:3 cm) node[VertexStyle, font=\scriptsize] (3) {$\lsb 1, 2, 3 \rsb$} ;
      \draw (65.1:3 cm) node[VertexStyle, font=\scriptsize] (1) {$\lsb 0, 1, 2 \rsb$} ;

      \draw [EdgeStyle] (0)--(1) ;
      \draw [EdgeStyle] (0)--(13) ;
      \draw [EdgeStyle] (0)--(5) ;

      \draw [EdgeStyle] (2)--(1) ;      
      \draw [EdgeStyle] (2)--(3) ;      
      \draw [EdgeStyle] (2)--(5) ;      

      \draw [EdgeStyle] (4)--(9) ;      
      \draw [EdgeStyle] (4)--(3) ;      
      \draw [EdgeStyle] (4)--(1) ;      

      \draw [EdgeStyle] (10)--(11) ;      
      \draw [EdgeStyle] (10)--(9) ;      
      \draw [EdgeStyle] (10)--(3) ;      

      \draw [EdgeStyle] (8)--(13) ;      
      \draw [EdgeStyle] (8)--(11) ;      
      \draw [EdgeStyle] (8)--(7) ;
      
      \draw [EdgeStyle] (6)--(13) ;      
      \draw [EdgeStyle] (6)--(7) ;      
      \draw [EdgeStyle] (6)--(5) ;      

      \draw [EdgeStyle] (12)--(11) ;      
      \draw [EdgeStyle] (12)--(7) ;      
      \draw [EdgeStyle] (12)--(9) ;      
    \end{tikzpicture}
  \end{center}
  \caption{$CH \lp \mathbb{Z}_7, \lsb \lsb 0, 1, 2 \rsb, \lb 0, 1, 6 \rsb, \lsb 0, 5, 6 \rsb \rsb \rp$}\label{nonLinear}
  \end{figure}

\begin{ex}
    Let $G=\mathbb{Z}_{7}$, and let $C_1 = \{0,1,2\},$ $C_2= \{0,1,6\},$ and $C_3=\{0,5,6\}$. The set $\mathcal{C}=\{C_1, C_2, C_3\}$ satisfies the \nameref{abc}, but the element $1$ is contained in two cells, while the element $2$ is only contained in one cell. The incidence graph of the resulting Cayley hypergraph is shown in Figure~\ref{nonLinear}.
\end{ex}

If we want to consider the class of graphs where the hypotheses of  Theorem~\ref{hyperSabidussi} are satisfied and an analogue of Proposition~\ref{cayleyHalved} is possible, we would need to restrict ourselves to hypergraphs where every pair of vertices is contained in either $0$ or $\lambda$ hyperedges. This is an unnatural constraint to place on hypergraphs unless $\lambda=1$, in which case the hypergraphs are linear. Linear hypergraphs are relevant to yet another hypergraph extension of Cayley graphs which we discuss in Section~\ref{secCages} and many of our later results will depend on the linearity of the hypergraph. 

Even though they were introduced in a very different context, the recent bipartite graphs of Evra, Feigon, Maurischat, and Parzanchevski~\cite{evra2023ramanujan} corresponds to a natural restriction of hypergraphs introduced by Shee~\cite{shee1990group}, Lee and Kwon~\cite{lee2013cayley}, and Jajcay and Jajcayov\'a~\cite{jajcayova2024generalizations}. Hypergraphs give another perspective through which we can view Cayley incidence graphs, and this perspective will be particularly useful in Sections~\ref{secCages} and~\ref{secAuth}.

\section{Designs and Difference Sets}\label{secDifference}

Let \( J \) be the all-ones matrix.
For a group $G$, $C \subseteq G$ and $g \in G$, define
\[
Cg := \{ sg : s \in C \},\quad C^{-1} := \{ s^{-1} : s \in C \}.
\]
A 2-\((v, k, \lambda)\) \textit{design} (or a \emph{$2$-design}) is a point-block incidence structure with $v$ points where every block contains $k$ points and every pair of points is contained in $\lambda$ blocks.
Each point belongs to $r = \lambda (v - 1) / (k - 1)$ blocks, and the total number of blocks is \(b = vr / k\).
If $N$ is the $v \times b$ incidence matrix, this definition is equivalent to the statement
\[ NN^T = \lp r - \lambda \rp I + \lambda J. \]
A 2-design with $\lambda = 1$ is often referred to as a \textit{Steiner system} \( S \lp 2, k, v \rp. \)
A \textit{symmetric 2-design} is a 2-\((v, k, \lambda)\) design with \(b = v,\) or equivalently, \(k = r\).
Symmetric 2-$(v, k, 1)$ designs are the point-line incidence structures of finite projective planes. Generally, constructions from finite geometry refer to lines instead of blocks, though for the purposes of this paper we will refer to them as blocks. For more background on designs, we refer the reader to the textbooks by Beth, Jungnickel, and Lenz~\cite{beth1999design} or Hughes and Piper~\cite{hughes1985design}.

Let $X = \BCay \lp G, \pi \rp$ be a Cayley incidence graph such that the underlying Cayley graph is complete, that is, $S \lp \pi \rp = G \setminus \lsb e \rsb.$
Letting $n$ be the order of $G$ and $k$ be the size of the cells of $\pi$, Proposition~\ref{cayleyHalved} tells us that $X$ is the point-block incidence graph of a 2-$\lp n, k, 1 \rp$ design.
Note that $G$ acts regularly on the points of this design.
Conversely, by Corollary~\ref{linearSabidussi}, the incidence graph of any 2-$(n, k, 1)$ design with an automorphism group acting regularly on its points is a Cayley incidence graph with a complete underlying Cayley graph.

When both halved graphs of a Cayley incidence graph $\BCay(G, \pi)$ are complete, we get the incidence graph of a projective plane with $G$ acting regularly on its points. The \nameref{Z7} of Example~\ref{Z7} is one such example.
It is known~\cite{beth1999design} that the existence of a projective plane (a symmetric 2-$(n, k, 1)$ design) which has $G$ as a point-regular automorphism group is equivalent to the existence of an \textit{$\lp n, k, 1 \rp$-difference set} of $G$.
This implies a correspondence between these difference sets and Cayley incidence graphs with complete halved graphs, and we further show this correspondence in a more direct way.

Note that the only known examples of finite projective planes with a group acting regularly on their points are Desarguesian planes.
It is a long standing open question~\cite{gill2016transitive} if there even exist non-Desarguesian finite projective planes with a point-transitive group action.
In 1959, Ostrom and Wagner~\cite{ostrom1959projective} proved that a finite projective plane with a group acting doubly transitively on its points is always Desarguesian.

Let $G$ be a group of order $n$.
A $k$-subset $D$ of $G$ is an $\lp n, k, \lambda \rp$-\textit{difference set} if every non-identity element of $G$ has exactly $\lambda$ representations as a difference $d \lp d' \rp^{-1}$ of two elements from $D$.
The term ``difference'' is more natural in additive notation. We, however, continue to use multiplicative notation since the group $G$ may be non-abelian.
The \textit{development} $\dev(D)$ of a difference set $D$ is the incidence structure with point set $G$ and block set $\{ Dg : g \in G \}$.
The development of an $(n, k, \lambda)$-difference set is a symmetric 2-$(n, k, \lambda)$ design.
More information on difference sets can be found in Beth, Jungnickel, and Lenz~\cite{beth1999design} or Jungnickel, Pott, and Smith~\cite{jungnickel200618difference}.

\begin{lem}\label{lm-diff-inv}
Let $G$ be a group, and let $D = \{d_1, \dots, d_k\}$ be a $k$-subset of $G$.
Then the following two conditions are equivalent:
\begin{enumerate}
\item for all $1 \leq i, j \leq k$ and $i \neq j$, the elements $d_id_j^{-1}$ are pairwise different;
\item for all $1 \leq i, j \leq k$ and $i \neq j$, the elements $d_i^{-1}d_j$ are pairwise different.
\end{enumerate}
\end{lem}
\proof
Suppose the first condition does not hold, i.e. $d_id_j^{-1} = d_sd_t^{-1}$ for some $i \neq j$, $s \neq t$ and $(i, j) \neq (s, t)$.
Note that $i = s$ would imply $j = t$, so we have $i \neq s$ and $j \neq t$.
Then $d_s^{-1}d_i = d_t^{-1}d_j$ for $s \neq i$, $t \neq j$ and $(s, i) \neq (t, j)$, thus the second condition also does not hold.
The same argument works in the other direction. \qed

If $D$ is a difference set, then so is $Dg$ for any $g \in G$, and, moreover, $\dev(D) = \dev(Dg)$. In particular, by possibly replacing $D$ with $Dd^{-1}$ for some $d \in D$, we may assume without loss of generality that a difference set $D$ contains $e$.

\begin{prop}\label{prop-diffset-cin}
Let $G$ be a group of order $n$, let $D$ be a $k$-subset of $G$ such that $e \in D$, let $\pi_{D^{-1}} = \{ dD^{-1} : d \in D \}$ and $\pi_D = \{ d^{-1}D : d \in D \}$.
Note that $\pi_{D^{-1}}$ and $\pi_D$ satisfy the \nameref{abc}.
The following conditions are equivalent:
\begin{enumerate}[label=(\arabic*)]
\item\label{itm-diffset-D} $D$ is an $(n, k, 1)$-difference set.
\item\label{itm-diffset-Dinv} $D^{-1}$ is an $(n, k, 1)$-difference set.
\item\label{itm-diffset-pi} Each non-identity element of $G$ belongs to precisely one cell of $\pi_{D^{-1}}$.
\item\label{itm-diffset-piinv} Each non-identity element of $G$ belongs to precisely one cell of $\pi_D$.
\end{enumerate}
If these conditions hold, the point-block incidence graphs of $\dev(D)$ and $\dev(D^{-1})$ and the Cayley incidence graphs $\BCay(G, \pi_{D^{-1}}) = (\bipartB \cup \bipartC, E)$ and $\BCay(G, \pi_D) = (\bipartB' \cup \bipartC', E')$ are all isomorphic.
\end{prop}

\proof
Conditions~\ref{itm-diffset-D} and~\ref{itm-diffset-pi} hold when all elements $d(d')^{-1}$ for $d \neq d'$, $d, d' \in D$, are pairwise different, and conditions~\ref{itm-diffset-Dinv} and~\ref{itm-diffset-piinv} hold when the same is true for elements $d^{-1}d'$.
Thus they are all equivalent due to Lemma~\ref{lm-diff-inv}.
The maps 
\[
\begin{cases}
g \mapsto D^{-1}g\\
Dh \mapsto h
\end{cases},
\quad
\begin{cases}
g \mapsto g^{-1}\\
Dh \mapsto h^{-1}D^{-1}
\end{cases},
\quad
\begin{cases}
g \mapsto g^{-1}D\\
Dh \mapsto h^{-1}
\end{cases}
\]
define graph isomorphisms from the point-block incidence graph of $\dev(D)$ to the other three graphs. \qed

Note that in Proposition~\ref{prop-diffset-cin} points of $\dev(D)$ correspond to the blocks of $\dev(D^{-1})$ and to $\bipartB$ and $\bipartC'$.

The correspondence between difference sets and Cayley incidence graphs is more straightforward when $G$ is abelian: if this is the case, $gD = Dg$ for any $g \in G$, so, continuing to use notation from Proposition~\ref{prop-diffset-cin}, the point-block incidence graph of $\dev(D)$ is simply equal to $\BCay(G, \pi_D)$. The \nameref{Z7} was constructed by $\BCay(\Z_7, \pi_D)$ for $D = \{0, 1, 3\}$.

Abel and Buratti~\cite{abel200618difference} state that the existence of a 2-$(n, k, 1)$ design with an automorphism group acting regularly on its points is equivalent to the existence of what they call an \textit{$(n, k, 1)$-partial difference family}.
This implies a correspondence between these partial difference families and Cayley incidence graphs with complete underlying Cayley graphs.
This can also be shown more directly in a similar manner to Proposition~\ref{prop-diffset-cin}, though we omit the details.

We conclude this section by showing how finite projective and affine spaces (Desarguesian projective and affine planes in case of dimension 2) can be constructed as Cayley incidence graphs.
Both of these constructions can be restated in terms of partial difference families.

For a prime power $q$, denote by $\FF_q$, $\FF_q^*$ and $\FF_q^+$ the finite field of size $q$, its multiplicative group and its additive group, respectively.

Let $n \geq 2$, let $q$ be a prime power, and let $\alpha$ be a generator of $\FF_{q^n}^*$.
Since $\alpha$ is a root of some primitive polynomial $f(x)$ of degree $n$ in $\FF_q$, every element of $\FF_{q^n}$ can be written as a sum, 
\[\sum_{i = 0}^{n - 1} c_i\alpha^i,\quad c_i \in \FF_q,\] 
thus giving a representation of $\FF_{q^n}$ as an $n$-dimensional vector space over $\FF_q$.
Denote by $\mathcal{V}_1$ the set of all 1-dimensional subspaces of this vector space, i.e. sets of the form $\langle a \rangle = \{ ka : k \in \FF_q \}$ for some $a \in \FF_{q^n} \setminus \{0\}$.
Let \( \qbino{n}{q} = |\mathcal{V}_1| = (q^n - 1) / (q - 1) \).

The following example for prime $q$ appears in Abel and Buratti~\cite{abel200618difference}.

\begin{ex}\label{exAGnq}
Let $G := (\FF_q^+)^n$ and $\pi := \mathcal{V}_1$.
Note that for any two different $\langle a \rangle, \langle b \rangle \in \pi$ we have $\langle a \rangle \cap \langle b \rangle = \{0\}$.
Also, each $\langle a \rangle \in \pi$ is an additive subgroup of $G$, so $\pi$ trivially satisfies the \nameref{abc}.
The Cayley incidence graph $\BCay(G, \pi) = (\bipartB \cup \bipartC, E)$ is in fact the point-line incidence graph of the affine space $\mathrm{AG}(n, q)$ with points $\bipartB$, lines $\bipartC$, and $\pi$ corresponding to the set of lines passing through the origin.
Note that this construction is a particular instance of Example~\ref{exSubgroups}.
\end{ex}

The next example for $n = 3$ corresponds to Singer difference sets~\cite{singerTheoremFiniteProjective1938}.

\begin{ex}\label{exPGnq}
Let $n \geq 3$.
Denote by $\mathcal{V}_2$ the set of all 2-dimensional subspaces of $\FF_{q^n}$ as a vector space over $\FF_q$.
Since $\FF_q^*$ is a subgroup of $\FF_{q^n}^*$, we have $\alpha^k \in \FF_q^*$ if and only if $k$ is divisible by $\abs{\FF_{q^n}^*} / \abs{\FF_q^*} = \qbino{n}{q}$.
In particular, $\langle \alpha^i \rangle \in \mathcal{V}_1$ are pairwise different for $0 \leq i < \qbino{n}{q}$, so $\mathcal{V}_1 = \{ \langle \alpha^i \rangle : 0 \leq i < \qbino{n}{q}\}$.
Now $\mathcal{V}_1$ may be identified with $G := \Z_{\qbino{n}{q}}$, and $\mathcal{V}_2$ may be identified with certain $(q + 1)$-subsets of $G$.

Let $\pi$ be the collection of $(q + 1)$-subsets of $G$ corresponding to 2-dimensional subspaces containing $\langle \alpha^0 \rangle$, i.e. sets of the form $\langle \alpha^0, \alpha^i \rangle = \lsb k_1 \alpha^0 + k_2 \alpha^i : k_1, k_2 \in \mathbb{F}_q \rsb \in \mathcal{V}_2$ for some $0 < i < \qbino{n}{q}$.
By construction, each $C \in \pi$ contains $0$.
For two non-equal $\langle \alpha^0, \alpha^i \rangle, \langle \alpha^0, \alpha^j \rangle \in \mathcal{V}_2$ we have $\langle \alpha^0, \alpha^i \rangle \cap \langle \alpha^0, \alpha^j \rangle = \langle \alpha^0 \rangle$, so two different cells of $\pi$ have only $0$ in common.
Moreover, $\pi$ satisfies the \nameref{abc}: if $C \in \pi$, $s \in C \setminus \{0\}$, then $C$ corresponds to $\langle \alpha^0, \alpha^s \rangle$, thus $C - s$ corresponds to $\alpha^{-s}\cdot \langle \alpha^0, \alpha^s \rangle = \langle \alpha^{-s}, \alpha^0 \rangle \in \mathcal{V}_2$, so $C - s \in \pi$.
The Cayley incidence graph $\BCay(G, \pi)$ is the incidence graph for the inclusion relation between $\mathcal{V}_1$ and $\mathcal{V}_2$, so it is the point-line incidence graph of the projective space $\mathrm{PG}(n - 1, q)$.

The collection $\pi$ described above can be calculated explicitly if the polynomial $f(x)$ is known; for example, to get \nameref{Z7} of Example~\ref{Z7} we can take $n = 3$, $q = 2$ and $f(x) = x^3 - x - 1$.
\end{ex}

\section{\texorpdfstring{$\bipartC$}{beta}-transitive Cayley Incidence Graphs}\label{secIncMin}

From Corollary~\ref{linearSabidussi}, we know that for a Cayley incidence graph $\BCay \lp G, \pi \rp = \lp \bipartB \cup \bipartC, E \rp$ the action of $G$ on $\bipartB$ is regular. However, this action need not be regular, or even transitive, on $\bipartC$. 

\begin{ex}The graphs in Example~\ref{Z3d} clearly have two orbits on $\bipartC$ under $G_1 \times G_2$.    
\end{ex}

For any Cayley incidence graph $\BCay \lp G, \pi \rp$ with vertex set $\bipartB \cup \bipartC$, we will show that the orbits under the action of $G$ on $\bipartC$ are classified by a relation on $\pi$. Then we will consider when the action is regular within each orbit.

Define $gC_i \approx hC_j$ if and only if $g C_i$ and $h C_j$ belong to the same orbit under $G$. This is an equivalence relation, and
\begin{align*}
    C_i \approx C_j &\iff \exists g \in G : gC_j = C_i\\
    & \iff \exists g \in C_i : gC_j = C_i\\
    & \iff \exists g \in C_i : C_j = g^{-1} C_i,
\end{align*}
so $C_i \approx C_j$ if and only if there is a $g\in C_i$ such that $g^{-1} C_i = C_j$. 

Since $gC_i \approx C_i$ for any $g\in G$, we see that $g C_i \approx h C_j$ if and only if $C_i \approx C_j$. We then restrict the relation $\approx$ to $\pi=\{C_1, \dots, C_{\ell} \}$, since this contains the same information as the full relation on $\bipartC$.

Denote by $[C]$ the equivalence class of $C \in \pi$ under the relation $\approx$, that is
\[
[C] = \{ s^{-1}C ~ : ~ s \in C \}.
\]

\begin{defn}
A Cayley incidence graph $\BCay(G, \pi)$ is called \textit{\(\bipartC\)-transitive} if for every $C_i, C_j \in \pi$ there is some $s\in C_i$ such that $s^{-1}C_i = C_j$. Equivalently, $\BCay \lp G, \pi \rp$ is $\bipartC$-transitive if $\lb C \rb = \pi$ for some $C \in \pi$.
\end{defn}

Note that this is not equivalent to $\BCay(G,\pi)$ having at most two orbits on the vertices. By our definition, $\beta$-transitivity refers to the action of $G$ on $\beta$, not the action of the full automorphism group. 

\begin{ex}
    The graph in the Example~\ref{Z7}~(\nameref{Z7}) is $\bipartC$-transitive. More generally, the Cayley incidence graphs coming from the difference set constructions in Proposition~\ref{prop-diffset-cin} are $\bipartC$-transitive.
\end{ex}

Let $X = \BCay \lp G, \pi \rp$ be a Cayley incidence graph and let $C \in \pi.$ If $X$ is $\bipartC$-transitive, then $\lb C \rb = \pi.$ More generally, $\lb C \rb$ satisfies the \nameref{abc}, so $\BCay \lp G, \lb C \rb \rp$ is an induced subgraph of $X$ which $G$ acts transitively on. This leads to the following result.

\begin{prop}\label{gammaOrbits}
Let $\BCay(G, \pi) = (\bipartB \cup \bipartC, E)$ be a Cayley incidence graph and let $\pi$ consist of $t$ equivalence classes: $\pi = [C_1] \sqcup \dots \sqcup [C_t]$.
The orbits of elements of $\bipartC$ under the action of $G$ are precisely
\[
\bipartC_i := \{gC_i : g \in G\}, 1 \leq i \leq t.
\]
\end{prop}

\begin{rem}
    This allows us to decompose a Cayley incidence graph into $\bipartC$-transitive Cayley incidence graphs. We cannot do the same if we weaken the linearity assumption: let $\mathcal{C} =  \lsb C_1, \ldots, C_{\ell} \rsb$ such that for all $1 \leq i \leq \ell$ the cell $C_i$ has $k$ elements and contains the identity. Suppose further that $\mathcal{C}$ satisfies the \nameref{abc} and that every non-identity element in $\bigcup_{i=1}^{\ell} C_i$ lies in the same number of cells of $\mathcal{C}$. We can still define $\bipartC$-transitive, but the equivalence class $\lb C_i \rb$ might not preserve the property that non-identity elements appear in the same number of cells.
\end{rem}

Now for any $gC \in \bipartC$ we will compute its \textit{setwise stabilizer}
\[
\Stab_G(gC) := \{ x \in G : x(gC) = gC \}.
\]
Since $\Stab_G(gC)= g\Stab_G(C)g^{-1}$, it suffices to compute $\Stab_G(C)$ for $C \in \pi$.

\begin{prop}\label{stabilizer}
Let $\BCay \lp G, \pi \rp$ be a Cayley incidence graph. For any $C \in \pi$ we have
\begin{equation*}
   \Stab_G(C) = C \cap C^{-1}.
\end{equation*}
\end{prop}
\proof
Let $S := C \cap C^{-1}$;
our goal is to prove that $\Stab_G(C) = S$.
For any $s \in S$, we have $s C \in \pi$ and $s \in sC \cap C$, thus $sC = C$, so $s \in \Stab_G(C)$.
Conversely, let $g\in \Stab_G(C)$, then $gC = C$, and since $e\in C$, it follows that $g, g^{-1}\in C$, thus $g \in S$.\qed

\begin{cor}Let $X := \BCay \lp G, \pi \rp = \lp \bipartB \cup \bipartC, E \rp$ be a Cayley incidence graph. The action of $G$ on $\bipartC$ is regular if and only if $X$ is $\bipartC$-transitive and $\Stab_G(C) = \{e\}$ for some $C \in \pi$.\qed
\end{cor}

When the action of $G$ on $\bipartC$ is regular, we will call it \emph{$\bipartC$-regular}.

\begin{rem}\label{rm:lleqk}
    A $\bipartC$-transitive graph $\BCay(G, \pi)$ always has $l \leq k$ because $\pi = [C]$ for any $C \in \pi$ and so  $\ell = |\pi| = |[C]| \leq |C| = k$. Equality holds exactly when $\Stab_G \lp C \rp = \{e\}$, so $\bipartC$-regular graphs are always regular as graphs.
\end{rem}

A \textit{Dedekind group} is a group in which every subgroup is normal. In particular, abelian groups are Dedekind groups.

\begin{cor}\label{cellSubgroup}Let $\BCay \lp G, \pi \rp$ be a Cayley incidence graph for some Dedekind group $G$. For any $C\in \pi$ either $\Stab_G(C) = C$ or $\Stab_G(C) = \{e\}$.
\end{cor}
\proof
Suppose that $\Stab_G(C) \neq C$, so by Proposition~\ref{stabilizer}, there exists some $s\in C$ with $s^{-1} \notin C$.
By the \nameref{abc}, $C' := s^{-1}C \in \pi$, and $s^{-1} \notin C$ implies $C \neq C'$.
Since $\Stab_G(C)$ is a normal subgroup of $G$, we have 
\[
\Stab_G(C') = s^{-1}\Stab_G(C)s = \Stab_G(C).
\]
But then $\Stab_G(C) \subseteq C \cap C' = \{e\}$.\qed

We conclude this section with a few more results on Cayley incidence graphs for Dedekind and abelian groups.

\begin{lem}\label{lemAbellgeqk}
Let $G$ be a Dedekind group, and let $\BCay(G, \pi)$ be an $(\ell, k)$-biregular $\bipartC$-transitive Cayley incidence graph.
Then $\ell = 1$ or $\ell = k$.
\end{lem}
\proof 
Fix any $C \in \pi$, then $\pi = [C] = \{s^{-1}C : s \in C\}$.
If $\ell \neq k$, then, by Remark~\ref{rm:lleqk}, $\ell < k$ and $\Stab_G \lp C \rp \neq \lsb e \rsb,$ so by Corollary~\ref{cellSubgroup}, $\Stab_G(C) = C.$ Thus $C$ is a subgroup, which implies $s^{-1}C = C$ for every $s \in C$. Since $\pi = \lb C \rb = \lsb C \rsb$, we must have $\ell =1$.\qed

\begin{lem}\label{girthSix}Let $X = \BCay \lp G, \pi \rp$ be an $(\ell, k)$-biregular $\bipartC$-transitive Cayley incidence graph with $G$ abelian. If $\ell \geq 3$, then $X$ has girth six.
\end{lem}

\proof By Lemma~\ref{girth6Plus}, $X$ has girth at least six. Let $C_1, C_2, C_3 \in \pi$ be pairwise distinct.
Since $\BCay(G, \pi)$ is $\bipartC$-transitive, there are some $s \in C_1$, $t \in C_2$ such that $C_2 = s^{-1}C_1$ and $C_3 = t^{-1}C_2$, so \[ sC_3 = st^{-1} C_2 = t^{-1} C_1 \] 
since $G$ is abelian. Thus
\[ e \lp C_1 \rp s \lp sC_3 \rp t^{-1} \lp C_3 \rp e \]
is a cycle of length six. \qed

Note that the same result does not hold if $G$ is non-abelian.

\begin{ex}\label{exNontransitive}
     Let $G$ be the non-abelian group of order 21. It is the semidirect product of $\Z_7$ and $\Z_3$ so we can write $G = \langle a,b\rangle$ where $a$ has order seven and $b$ has order three. Let $\pi = \{C_1, C_2, C_3\}$ with
     \begin{align*}
         C_1 & = \{e, b,b^{-1}\},\\
         C_2 & = \{e, ab,(ab)^{-1}\}\\
         C_3 & = \{e, a^2b, (a^2b)^{-1}\}
     \end{align*}
     Then $\pi$ satisfies the \nameref{abc} and $k= \ell = 3$. The Cayley incidence graph is regular but not a Cayley graph and it has girth eight.
\end{ex}

\section{Subdivision Graphs and Cages}\label{secCages}

Whenever we have a bipartite biregular graph where one bipart has valency two, we can consider it as the subdivision graph of a regular graph. We used this in Example~\ref{exDiscrete} to show that the subdivision graph of any Cayley graph is a Cayley incidence graph.

However, if $\pi$ has two cells, then the Cayley incidence graph is the subdivision graph of the other bipart and not the underlying Cayley graph. Thus we wish to characterize the Cayley incidence graphs \( \BCay \lp G, \lsb C_1, C_2 \rsb \rp.\)

\begin{prop}\label{l2}
Let \( X = \BCay \lp G, \lsb C_1, C_2 \rsb \rp \) be a Cayley incidence graph. Then there exist subgroups \( S_1, S_2\) with trivial intersection such that either
\begin{enumerate}[label=(\arabic*)]
\item \label{l2one} \( C_1 = S_1 \) and \( C_2 = S_2\), or
\item \label{l2two} there exists some non-identity group element \( x \) such that $S_2 = x^{-1} S_1 x$, $C_1= S_1 \cup S_1 x$ and $C_2 = S_2 \cup S_2 x^{-1}.$
\end{enumerate}
\end{prop}
\proof If $X$ is not $\bipartC$-transitive, then clearly~\ref{l2one} holds. Thus we may assume there exists some $x \in C_1$ such that $x^{-1}C_1=C_2$. Let $S_1 = \Stab_G(C_1)$ and observe that $S_1x\subseteq C_1$. To see that these are all the elements of $C_1$, we suppose that there is some other element $y$ such that $y^{-1}C_1 = C_2$. Then $yx^{-1} \in S_1$ since $yx^{-1}C_1 = yC_2 = C_1$, and thus $y\in S_1x$. Similarly, letting $S_2 = \Stab_G(C_2)$, we get $C_2 = S_2 \cup S_2x^{-1}$, and by the property of stabilizers, $C_2 = x^{-1}C_1$ implies $S_2 = x^{-1}S_1x$.\qed

For abelian groups, this gives us a very restricted form.

\begin{cor}Let $G$ be an abelian group of order $n$ and let $X = \BCay \lp G, \lsb C_1, C_2 \rsb \rp$ be a connected $\lp 2, k \rp$-biregular Cayley incidence graph. Then either $k=2$ and $X$ is a cycle, or there exists some $m$ such that $n = mk$ and $X$ is the subdivision graph of the complete bipartite graph \( K_{m,m} \).
\end{cor}

\proof By Corollary~\ref{cellSubgroup}, the only way Case~\ref{l2two} of Proposition~\ref{l2} can hold is if $S_1 = S_2 = \{e\}$, in which case we have $\BCay \lp G, \lsb \lsb e, x \rsb, \lsb e, x^{-1} \rsb \rsb \rp$, which is a cycle.

Otherwise, $C_1$ and $C_2$ are two disjoint subgroups of order $k$. Then the halved graph $H_\bipartC$ is the bipartite graph with biparts given by the cosets $G / C_1$ and $G / C_2$. Since the elements of $C_1 \cup C_2$ generate $G$, we must have $G \cong C_1 \times C_2.$ Thus $$m:= \abs{G/C_1} = \abs{C_1} = \abs{C_2} = \abs{G/C_2}$$ is both the valency of vertices  in $H_{\bipartC}$ and the sizes of the biparts. It follows that $H_{\bipartC}$ is the complete bipartite graph $K_{m,m}$ and $X$ is its subdivision graph. \qed

For non-abelian groups, the examples are less trivial.
Given two subgroups $S_1, S_2$, the \textit{rank-two coset geometry} can be viewed as a bipartite graph with biparts formed from the cosets of $G/S_1$ and $G/S_2$. The classification of coset geometries was spearheaded by Buekenhout~\cite{buekenhout2006geometry, buekenhout1995finite}. 

Many results in this area are concerned with classifying rank-two coset geometries for a particular group. The incidence graph of a rank-two coset geometry is always bipartite and edge-transitive. De Saedeleer, Leemans, Mixer, and Pisanki~\cite{saedeleer2014core} took the opposite approach and studied regular, edge-transitive bipartite graphs of valency three and four to see which graphs were the incidence graphs of rank-two coset geometries.

Case~\ref{l2one} in Proposition~\ref{l2} is the subdivision graph of the incidence graph of a rank-two coset geometry with $\abs{S_1} = \abs{S_2}$ and $S_1$ and $S_2$ intersecting only in the identity.

\begin{ex}\label{cosetGeo}Let $G$ be the unique non-abelian group of order 21, that is $G \cong \Z_7\rtimes\Z_3$. Let $a$ and $b$ be generators of order three and seven, respectively and define $C_1:=\langle a\rangle$ and $C_2:=\langle ab\rangle$. Then $C_1$ and $C_2$ are subgroups of order three with trivial intersection and the graph $\BCay(G,\{C_1,C_2\})$ is a $(2,3)$-biregular graph of girth twelve. The halved graph $H_{\bipartC}$ is isomorphic to the Heawood graph, shown in Figure~\ref{heawood}. \end{ex}

We also have non-trivial examples of Case~\ref{l2two} in Proposition~\ref{l2}.

\begin{ex}\label{exl2S4}
Consider $G= S_4$, and let $C_1=\{e, (12), (134), (1342)\}$ and $C_2=\{e, (24), (431), (1243)\}$. This gives an example of Case~\ref{l2two} where $x_1 = (134)$ and $S_1=\{e, (12)\}, S_2=\{e, (24)\}$.
\end{ex}

A $\lp k, g \rp$ \textit{cage} is the graph with the minimum number of vertices over all $k$-regular graphs of girth $g.$ We can similarly define an $\lp \ell, k, g \rp$-biregular cage as the bipartite graph with minimum number of vertices over all $\lp \ell, k \rp$-biregular graphs of girth $g.$ This is equivalent to viewing them as hypergraph cages as done by Ellis and Linial~\cite{ellis2014regular} or Erskine and Tuite~\cite{erskine2023small}. More information about regular cages can be found in the survey of Exoo~\cite{exoo2012dynamic}, but one major idea is that Cayley graphs can be used to construct graphs of given valency and large girth. 

We could use Proposition~\ref{l2} to develop a strategy to search for $\lp 2, 3, g \rp$ or $\lp 2, 4, g \rp$-cages using Cayley incidence graphs.

In Case~\ref{l2one}, the work of De Saedeleer et al.~\cite{saedeleer2014core} gives coset geometries with $3$- and $4$-regular incidence graphs. Using the subgroups from these coset geometries as the cells of our Cayley incidence graphs gives us $\lp 2, 3 \rp$ and $\lp 2, 4 \rp$--biregular graphs. Thus, if we begin by looking at cubic and quartic bipartite edge-transitive graphs with large girth, we can construct $\lp 2, 3 \rp$ and $\lp 2, 4 \rp$-biregular graphs with large girth.

In Case~\ref{l2two}, we can build the cells of the Cayley incidence graph in a deliberate way from the groups. By Lemma~\ref{girthSix}, we know we can restrict our search to non-abelian groups, and then from Proposition~\ref{l2} we look for nontrivial subgroups of order $\frac{k}{2}$ in the same conjugacy class. This cannot happen when $k=3$, and when $k=4$ the desired subgroups have order $2$.

This would be a computationally viable, but mathematically inefficient, way to search for $\lp 2, k \rp$-biregular Cayley incidence graphs of large girth. One of the steps involved is to look at $k$-regular graphs with large girth and impose additional restrictions for which graphs we can construct as Cayley incidence graphs. More generally, for any $k$-regular graph with large girth, we can construct a $\lp 2, k \rp$-biregular bipartite graph with twice the girth simply by subdividing every edge once. However, a similar strategy for $\lp 3, k \rp$-biregular Cayley incidence graphs should lead to non-trivial constructions.

This is the general approach used by Erskine and Tuite~\cite{erskine2023small} and Ellis and Linial~\cite{ellis2014regular} in considering hypergraph cages. However, rather than using Cayley incidence graphs, they used another hypergraph generalization of Cayley graphs due to Buratti~\cite{buratti1994cayley}.

\begin{defn}[$t$-Cayley hypergraph]\label{tCayley}
Let \( G \) be a group, let \( S \) be a subset of \( G \) that does not contain the identity, and let $t$ be an integer between two and the maximum order of elements of $S$. Then the $t$-\textit{Cayley hypergraph} has vertex set \( G \) and hyperedge set
\[ \lsb g s^i : 0 \leq i \leq t-1, g \in G, s \in S \rsb.\]    
\end{defn}

For a set \( S \subseteq G \setminus \lsb e \rsb, \) we can define a collection of subsets
\[ \mathcal{C} := \lsb \lsb e, s, s^2, \ldots, s^{t-1} \rsb : s \in S \rsb, \]
from which it is clear that $t$-Cayley hypergraphs are a subclass of the more general class of group hypergraphs in Definition~\ref{groupHypergraph}. If $t$ is at most minimum order of elements of $S$, then the $t$-Cayley hypergraph is $t$-uniform.

Buratti~\cite{buratti1994cayley} verified that a $t$-Cayley hypergraph has a group of automorphisms which act regularly on the vertices, but also showed that the Fano plane is a 3-uniform hypergraph which has a group of automorphisms acting regularly on the vertices, even though it is not a 3-Cayley hypergraph. Notably, the \nameref{Z7} is the Cayley incidence graph of Example~\ref{Z7} and a $\lp 3, 6 \rp$-cage, or $\lp 3, 3, 6 \rp$-biregular cage.

Ellis and Linial~\cite{ellis2014regular} independently defined a construction of random Cayley hypergraphs that is a subconstruction of $t$-Cayley hypergraphs. In the notation of Buratti~\cite{buratti1994cayley}, Ellis and Linial~\cite{ellis2014regular} let \( G \) be the symmetric group and let \( S \) be a set chosen uniformly at random from permutations with a particular structure. They proved that with high probability, the ensuing $t$-Cayley hypergraph is a linear, uniform, regular hypergraph with large girth.

A related, but more explicit, construction using $t$-Cayley hypergraphs was given recently by Erskine and Tuite~\cite{erskine2023small}, who used Definition~\ref{tCayley} to obtain computational results improving the upper bound on the size of hypergraph cages. We propose that Cayley incidence graphs are precisely the right class of graphs to search for biregular cages. The restriction of Definition~\ref{groupHypergraph} to $t$-Cayley graphs excludes important graphs, like the \nameref{Z7}, but by Corollary~\ref{linearSabidussi} we see that the Cayley incidence graphs are precisely the restriction of Definition~\ref{groupHypergraph} that might give rise to biregular cages.

\section{Automorphisms}\label{secAuth}

We have already seen how a group $G$ acts on $\BCay(G,\pi)$. In this section we will further investigate the full automorphism group of a Cayley incidence graph by considering how group automorphisms may yield additional automorphisms of the graph.

One notable use of Cayley graphs was by Frucht~\cite{frucht1939herstellung} to prove that every group is the automorphism group of a finite undirected graph. The idea of the proof and subsequent strengthenings was to begin with a Cayley graph and replace edges with certain subgraphs to disrupt the symmetries and obtain the desired group. This is related to the problem of \textit{graphical regular representations}, i.e. Cayley graphs with $\Aut \lp \Cay \lp G, S \rp \rp \cong G$. 

Godsil~\cite{godsil1981grrs} proved that with the exception of two infinite families and 13 simple groups with order at most 32, every group admits a graphical regular representation. Babai, Godsil, Imrich, and Lovasz~\cite{babai1982automorphism} conjectured that if $(G_n)$ is a family of groups with order $n$ tending to infinity that is neither abelian with exponent greater than two nor generalized dicyclic, then almost all Cayley graphs are graphical regular representations. This conjecture was proven recently by Xia and Zheng~\cite{xia2023asymptotic}.

When considering the automorphism group of a Cayley incidence graph $X = \BCay \lp G, \pi \rp = \lp \bipartB \cup \bipartC, E \rp$, there is some choice to be made. Viewing $X$ as a bipartite graph, an automorphism maps $\bipartB \cup \bipartC$ to $\bipartB \cup \bipartC$ and preserves the adjacency relation between $\bipartB$ and $\bipartC$. However, viewed as a hypergraph, an automorphism maps $\bipartB$ to itself and preserves the hypergedge relations. In particular, if a Cayley incidence graph is regular, we might have a graph automorphism mapping $\bipartB$ to $\bipartC$ which is not a hypergraph automorphism.

For any Cayley incidence graph we have $\Aut \lp CH \lp G, \pi \rp \rp \subseteq \Aut \lp \BCay \lp G, \pi \rp \rp$, and so we will begin by considering the hypergraph automorphism group. 

Regular representations of hypergraphs were considered by Bayat, Alaeiyan, and Firouzian~\cite{bayat2019normality} for $t$-Cayley hypergraphs, and by Jajcay~\cite{jajcay2002representing} for group hypergraphs. We are specifically interested in Cayley incidence graphs and a result used by Godsil~\cite{Godsil1981} to prove that almost all Cayley graphs of certain groups are graphical regular representations.

Given a graph $X$, and a group $G$ acting transitively on the vertices of $X$, Godsil~\cite{Godsil1981} determined the normalizer of $G$ in $\Aut(X)$. We similarly determine the normalizer of the group in the automorphism group of a Cayley incidence graph.

We say that a group automorphism $\varphi\in \Aut(G)$, \emph{induces a permutation of $\pi= \{C_1, \dots, C_\ell\}$}, if for each $1 \leq i \leq \ell$ there exists some $j\in \{1, \dots, \ell\}$ such that $\varphi(C_i) = C_j$. Let
\begin{equation*}
    \Aut(G, \pi) := \{\varphi \in \Aut(G): \varphi \text{ induces a permutation of }\pi\}.
\end{equation*}

For a group $G$ and a subgroup $G'$, we will denote by $N_G(G')$ the normalizer of $G'$ in $G$. Let $AH = \Aut \lp CH \lp G, \pi \rp \rp.$ 
For $g \in G$, let $\mathcal{L}_g \in \Aut(G)$ be the left translation $\mathcal{L}_g: G\rightarrow G: x\mapsto gx$.
The set $\{\mathcal{L}_g : g \in G\}$ forms a subgroup of $\Aut(G)$ isomorphic to $G$. We will abuse the notation and refer to this subgroup as $G$.

\begin{lem}\label{lem-fhom}
    Let $G$ be a group, and let $f: G\rightarrow G$ be a bijective function with $f(e)=e$, such that for all $g\in G$,
    \begin{equation*}
       f\circ  \mathcal{L}_g \circ f^{-1}=\mathcal{L}_h
    \end{equation*} 
    for some $h \in G$.
    Then $f$ is a group automorphism.
\end{lem}
\proof
We see that $f(g) = f(gf^{-1}(e)) = h$. Thus, we have
\begin{equation*}
    f \circ \mathcal{L}_g \circ f^{-1} =  \mathcal{L}_{f(g)},
\end{equation*}
and therefore,
\begin{equation*}
  \mathcal{L}_{f(gh)} =   f \circ \mathcal{L}_{gh} \circ f^{-1}   =  f \circ \mathcal{L}_{g} \circ f^{-1} \circ  f \circ \mathcal{L}_{h} \circ f^{-1}=   \mathcal{L}_{f(g)} \circ  \mathcal{L}_{f(h)}.
\end{equation*}
This implies $\mathcal{L}_{f(gh)}(e) =   \mathcal{L}_{f(g)} \circ  \mathcal{L}_{f(h)}(e)$, so we have $f(gh)=f(g)f(h)$. Then, since $f(e)=e$, $f$ is a group automorphism.
\qed

\begin{thm}\label{auts_hyp}
    Let $\BCay(G, \pi)$ be a Cayley incidence graph. Then 
    \begin{equation*}
    N_{AH}(G) \cap \{\varphi \in AH : \varphi(e)=e \}  \cong \Aut(G, \pi).
    \end{equation*}
\end{thm}
\proof We recall that the vertex set of $CH \lp G, \pi \rp$ is given by $G$. Hence, hypergraph automorphisms $f\in N_{AH}(G)$ are in particular bijections $f: G\rightarrow G$, and $G$ is a subgroup of $\Aut\lp CH \lp G, \pi \rp\rp$.

Let $f: G\rightarrow G$ be a hypergraph automorphism with $f \in N_{AH} \lp G \rp$ and $f(e)=e$. Since $f\in N_{AH}(G)$, for every $g\in G$ there exists some $h \in G$ such that
\begin{equation*}
    f \circ \mathcal{L}_g \circ f^{-1} =  \mathcal{L}_h.
\end{equation*}
Thus, by Lemma~\ref{lem-fhom}, $f$ is a group automorphism.
Since $f(e)=e$, we also see that $f$ must induce a permutation of $\pi$, by the fact that it is a hypergraph automorphism. Hence $f\in \Aut(G, \pi)$.

We now show that any $\varphi\in \Aut(G, \pi)$ defines a hypergraph automorphism and that $\varphi \in N_{AH}(G)$. Let $gC_i$ be a hyperedge. Then $\varphi(gC_i) = \varphi(g) \varphi(C_i)$  and hence the result is a hyperedge. Since $\Aut(G, \pi)$ is a group, we see that the preimage of any hyperedge is also a hyperedge. Thus elements of $\varphi\in \Aut(G,\pi)$ define automorphisms of $\HCay(G,\pi)$. 

For any $h\in G$, we have
\begin{equation*}
    \varphi \circ \mathcal{L}_g\circ \varphi^{-1}(h) = \varphi(g)h,
\end{equation*}
and hence $\varphi \circ \mathcal{L}_g\circ \varphi^{-1} = \mathcal{L}_{\varphi(g)}$. Thus $\varphi \in N_{AH}(G).$\qed

\begin{rem}
    Theorem \ref{auts_hyp} has a straightforward generalization to the automorphism groups of group hypergraphs. The proof is identical.
\end{rem}

\begin{cor}
Let $\BCay(G, \pi)$ be a Cayley incidence graph. Then
$$N_{AH}(G) \cong G \rtimes \Aut(G, \pi).$$
\end{cor}

\proof
    Let $f\in N_{AH}(G)$. Suppose that $f(e) = g$. Then $\mathcal{L}_{g^{-1}} \circ f(e)=e$, and thus $\mathcal{L}_{g^{-1}} \circ f \in \Aut(G, \pi)$ by the previous proposition. It is now clear that any element of $N_{AH}(G)$ can be written as a product of $g\in G$ and $f\in \Aut(G, \pi)$, and it is clear that $G\cap\Aut(G, \pi)=\{e\}$. Thus $N_{AH}(G)$ is the stated semi-direct product.
\qed

Now let us consider the existence of automorphisms which swap the two sides of $\BCay(G, \pi)$. In terms of hypergraphs, we can think of these as isomorphisms between a hypergraph and its dual. However, these are truly automorphisms of the bipartite graph, and it will be convenient to therefore view them as graph automorphisms with $A = \Aut \lp \BCay \lp G, \pi \rp \rp.$ We note that one can only have an automorphism $f\in N_{A}(G)$ which swaps $\bipartB$ and $\bipartC$ if the graph is $\bipartC$-regular, so we restrict our attention to this case.

To understand when an automorphism can swap two sides of a Cayley incidence graph, we introduce the following construction.

\begin{prop}\label{Dual}
    Let $X=\BCay(G, \pi)$ be $\bipartC$-regular with $\pi=\{C_1, C_2, \dots, C_k\}$, $C_1= \{g_1, g_2, \dots, g_k\}$, $g_1 = e$ and $C_i = g_i^{-1}C_1$. Let $\pi^* = \lsb C_1^{-1}, g_2C_1^{-1}, \ldots, g_k C_1^{-1} \rsb$. Then $\BCay \lp G, \pi^* \rp$ is a well-defined Cayley incidence graph isomorphic to $X$.
\end{prop}
\proof The cells $C_1^{*}, \ldots C_k^{*}$ of $\pi^*$ are given by
\begin{align*}
    C_1^{*} &=\{e, g_2^{-1}, \dots, g_k^{-1} \}\\
    C_2^{*} &=\{g_2, e, \dots, g_2 g_k^{-1} \}\\
&\vdots\\
    C_k^{*} &=\{g_k, g_kg_2^{-1}, \dots, g_k g_{k-1}^{-1}, e \}.
\end{align*}
It is straightforward to check that $\pi^{*}$ satisfies the \nameref{abc}. We additionally need to check that the cells $C_i^{*}=\{g_i g_j^{-1} ~: ~ g_j\in C_1\}$ only intersect in the identity. Since $C_i= g_i^{-1}C_1 = \{g_i^{-1}g_j ~: ~ g_j\in C_1\}$, and since two cells $C_i, C_{a}\in\pi$ only intersect in the identity, we see that the elements $g_i^{-1}g_j$ are all pairwise distinct when $i\neq j$. Hence, by Lemma \ref{lm-diff-inv}, the elements $g_i g_j^{-1}$ are all distinct. Thus two cells $C_i^{*}, C_{a}^{*}$ only intersect in the identity, and hence $\BCay(G,\pi^*)$ is well defined.

Let $\bipartB \cup \bipartC$ and $\bipartB^* \cup \bipartC^*$ be the vertex sets of $\BCay(G, \pi)$ and $\BCay(G, \pi^{*})$, respectively.
Let $\psi$ be the map $\bipartB \cup \bipartC \rightarrow \bipartB^* \cup \bipartC^*$ defined by  
  \begin{equation*}
        \begin{cases}
            g\mapsto gC_1^{*} \text{ for } g\in \bipartB\\
            gC_1\mapsto g \text{ for } gC_1\in \bipartC,
        \end{cases}
    \end{equation*}
Note that $h \in gC_1$ if and only if $\psi(gC_1) = g \in hC_1^* = \psi(h)$, thus $\psi$ is a graph homomorphism.
Since $(\pi^{*})^{*} = \pi$, it is easy to see that 
\begin{equation*}
    \varphi: G\rightarrow G: 
    \begin{cases}
        g\mapsto gC_1 \text{ for } g\in \bipartB^*\\
        gC_1^{*} \mapsto g \text{ for } g\in \bipartC^*,
    \end{cases}
\end{equation*}
defines a two-sided inverse of $\psi$. 
\qed

\begin{cor}
    Let $X=\BCay(G, \pi)$ be a $\bipartC$-regular Cayley incidence graph. Then both halved graphs of $X$ are Cayley graphs and if $S(\pi) = S(\pi^{*})$, the two halved graphs are isomorphic. \qed
\end{cor}

Note that for an abelian group $S(\pi) = S(\pi^*)$ always holds.

We will need the following standard result about bipartite graphs.

\begin{prop}\label{bipartite_auts}
Let $X$ be any connected bipartite graph with biparts $\bipartB$ and $\bipartC$. There is a well defined automorphism
$$\sigma: \Aut(X)\rightarrow \Sym{(\{\bipartB,\bipartC\})},$$
taking any graph automorphism to the permutation which it induces on $\{\bipartB, \bipartC\}$. Let $G_0=\ker(\sigma),$ i.e.:
\begin{equation*}
    G_0 = \{\varphi\in \Aut(X) : \varphi(\bipartB) =\bipartB,  \varphi(\bipartC) =\bipartC \}.
\end{equation*}
Then either $G_0 = \Aut(X)$ or $G_0$ has index 2 in $\Aut(X).$
\end{prop}
\proof
Let $\varphi$ be any automorphism of $X$. Since $\varphi$ preserves distances, we see that for any vertices $\vtxa, \vtxb$ at an even distance, the vertices $\varphi(\vtxa)$ and $\varphi(\vtxb)$ are at an even distance. Since $X$ is connected, two vertices are at an even distance if and only if they belong to the same set in the bipartition, and thus we see that $\varphi \in G_0$ or $\varphi(\bipartC) = \bipartB$ and  $\varphi(\bipartB) = \bipartC$.
Hence, the group homomorphism $\sigma: \Aut(X)\rightarrow \Sym{(\{\bipartB,\bipartC\})}$ is well defined, and since $\ker(\sigma)= G_0$, the result follows. 
\qed

\begin{thm}\label{Normalizer-duals}
    Let $\BCay(G, \pi) = \lp \bipartB \cup \bipartC, E \rp$ be a connected $\bipartC$-regular Cayley incidence graph with $\pi=\{C_1, C_2, \dots, C_k\}$,  $C_1= \{e, g_2, \ldots, g_k\}$, $g_1 = e$ and $C_i = g_i^{-1}C_1$. Let $\pi^* = \lsb C_1^*, \dots, C_k^* \rsb$, where $C_i^* = g_iC_1^{-1}$.
    Let
    \[
    \sigma: A\rightarrow \Sym{(\{\bipartB,\bipartC\})}
    \]
 send any graph automorphism to the permutation which it induces on $\{\bipartB, \bipartC\}.$
If there exists a group isomorphism $\varphi: G\rightarrow G$ such that $\varphi(C_1) = C_j^{*}$ for some $j$, then
    \begin{equation*}
     N_{A}(G) \cong \langle N_{AH}(G) , \varphi \rangle,
    \end{equation*}
    and $\ker(\sigma) = N_{AH}(G)$. Otherwise,
    \begin{equation*}
     N_{A}(G) \cong N_{AH}(G).
    \end{equation*}
\end{thm}
\proof We begin by showing that given a group automorphism $\varphi: G\rightarrow G$ with $\varphi(C_1) = C_{j}^{*}$, we can define an automorphism $f$ of the graph. Considering $C_1$ and $C_1^*$ as sets of group elements, we know there is some $g_{\varphi}$ such that $\varphi \lp C_1 \rp = g_{\varphi} C_1^*.$ We can thus define $f_{\varphi}: \bipartB \cup \bipartC \rightarrow \bipartB \cup \bipartC$ by
\begin{equation*}
\begin{cases}
g\mapsto \varphi(g)C_1 & \text{for }g\in \bipartB\\
g C_1 \mapsto \varphi(g) g_{\varphi} & \text{for }gC_1\in \bipartC\\
\end{cases}
\end{equation*}
We remark that
\begin{equation*}
f_{\varphi} \circ \mathcal{L}_{g} = \mathcal{L}_{\varphi(g)} \circ f_{\varphi},
\end{equation*}
so if the resulting map $f_{\varphi}$ is in  $A$,  then we have $f_\varphi \in N_{A}(G)$. Note that $h \in gC_1$ is equivalent to $\varphi(g^{-1}h) \in \varphi(C_1) = g_\varphi C_1^*$, which is in turn equivalent to $f_\varphi(gC_1) = \varphi(g)g_\varphi \in \varphi(h)C_1 = f_\varphi(h)$, so $f_\varphi$ is a graph automorphism.

Consider the map $\sigma: A\rightarrow \Sym{(\{\bipartB,\bipartC\})}$.
Since $N_A(G) \cap \ker(\sigma) = N_{AH}(G)$, by Proposition~\ref{bipartite_auts} we see that either $N_A(G)$ is a two-extension of $N_{AH}(G)$ or $N_{AH}(G) = N_A(G)$. We are in the first case if there exists any element $f\in N_A(G) \setminus N_{AH}(G)$, and in the second case otherwise. If we are in the second case, we are done, so we may suppose that we are given an automorphism $f \in N_A \lp G \rp$ that swaps the biparts.

Since $\BCay(G, \pi)$ is $\bipartC$-transitive, we can identify $\bipartC$ with $G$ through the bijection $g\mapsto g C_1$. Then, the restriction $f_{\vert \bipartB}$ can be viewed as a map $G\rightarrow G$, and we note that the restriction of the inverse $(f^{-1})\vert_{\bipartC}$ is the inverse of $f_{\vert \bipartB}$. After applying a translation $\mathcal{L}_g$, we may assume that $f(e)=e$. Viewed as a mapping of the vertices of the graph, this means that $f(e) = C_1$. Then, since $f \in N_A(G)$, we have $f_{\vert \bipartB} \circ \mathcal{L}_g \circ f_{\vert \bipartB}^{-1} = \mathcal{L}_h$, and thus, by Lemma~\ref{lem-fhom}, $f_{\vert \gamma}$ is a group automorphism. Since $f$ maps edges to edges, we see that any edge $e\sim C_i$ maps to an edge $f(e) \sim f(C_i)$. Note that this means that the neighbourhood of $C_1$, which consists of elements $\{e, g_2, \cdots, g_k\}\subseteq \bipartB$, must be mapped to neighbours of $f(C_1)$. Since $f(e) = C_1$, we see that we must map the vertex $C_1$ to exactly one of the vertices $e, g_2, \dots, g_k$. Thus $\{e, g_2, \dots, g_k\}$ must be mapped under $f$ to neighbours of $g_j$ for some $j$. The neighbours of $g_j$ are precisely
\begin{equation*}
g_{j}g_i^{-1}C_1,\text{ where } i\in \{1, \dots, k\}.
\end{equation*}
Using the identification $g C_1 \leftrightarrow g$, we see that $$f_{\vert \gamma}(\{e, g_2, \dots, g_k\}) = \{g_j g_{i}^{-1} : j\in \{1, \dots, k\}\},$$ so considering $C_1$ as a set of group elements, we have $f_{\vert \gamma}(C_1) = C_j^{*}$.\qed

\begin{ex}
    We consider the quaternion group $Q_8$, with $$\pi= \{\{1, i, -j\}, \{1, -i, k\}, \{1, j, -k\}\},$$ and we let $\varphi:Q_8\rightarrow Q_8$ be given by $i\mapsto j$, $j\mapsto i$, $k\mapsto -k$. Then $\pi'$ is given by
    \begin{equation*}
        \pi' = \{\{1, -i, j\}, \{1, i, k\}, \{1, -j, k\}\}
    \end{equation*}
    This induces an automorphism of $\BCay(Q_8, \pi)$ which swaps the two sides.  
\end{ex}

\section{Bipartite Cayley Graphs and Bi-Cayley Graphs}\label{secBiBi}

We constructed the \nameref{Z7} in Example~\ref{Z7} as a Cayley incidence graph, but we can also construct it as a regular Cayley graph. Consider the dihedral group of order $14$
\[D_7= \langle a,b: a^n=b^2=e, bab = a^{-1}\rangle,\] 
and let $T=\{b,ab,a^3b\}$. Then $\Cay \lp D_7, T \rp$ is isomorphic to $$\BCay \lp \mathbb{Z}_7, \lsb \lsb 0, 1, 3 \rsb, \lsb 0, 4, 5 \rsb, \lsb 0, 2, 6 \rsb \rsb \rp.$$
More generally, it is natural to ask when a regular Cayley incidence graph is a bipartite Cayley graph, or vice-versa.

We will need the following fact about bipartite Cayley graphs.

\begin{prop}\label{cayleyBipartite}
    A Cayley graph $\Cay \lp G, S \rp$ is bipartite if and only if there exists a homomorphism $\varphi: G \to \mathbb{Z}_2$ such that $\varphi \lp S \rp \subseteq \lsb 1 \rsb.$
\end{prop}
\proof
    Given such a homomorphism, it is clear that $\varphi^{-1}(0)$ and $\varphi^{-1}(1)$ define a bipartition of $V(\Cay(G,S))$.
    
    Conversely, suppose that we are given a Cayley graph which is bipartite. We consider the homomorphism $\varphi: G\rightarrow \Z_2$ which is the homomorphism $\sigma$ as in Proposition \ref{bipartite_auts}, restricted to $G$. It is surjective since $G$ acts transitively on $\Cay(G,S)$. By the way that elements of $S$ act on $\Cay(G,S)$, we also see that all elements of $S$ must swap the sides, and hence we have $\varphi(S)\subseteq \{1\}$.\qed

We begin by establishing that any bipartite Cayley graph with girth at least six is a Cayley incidence graph.

\begin{thm}\label{cayleyBiCayley}If $X = \Cay \lp G, S_1 \rp$ is bipartite with girth at least six, then $X \cong \BCay \lp G_0, \pi \rp$ for some subgroup $G_0 \leq G$ of index two and some $\pi$ with $S \lp \pi \rp = S_1$.
\end{thm}

\proof Since $X$ has no four-cycles, we know that for all $s_i \neq s_j \in S_1$ the products $s_i s_j$ must be distinct.

 Let $S_2= \{ s_i s_j \in G : 1 \leq i,j \leq \ell \}$, and define $\pi = \{C_1, \dots, C_\ell\}$ by setting $C_i=\{ s_is_j : 1 \leq j \leq \ell\}$. So 
\begin{align*}
C_1 &=\{s_1^{2},\dots, s_1s_\ell\},\\
C_2 &=\{s_2s_1, s_2^{2},\dots, s_2s_\ell\},\\
 & ~ ~ \vdots\\
C_\ell &= \{s_\ell s_1,\dots, s_\ell^{2}\}.
\end{align*}
Note that each non-identity element in $S_2$ appears exactly once, and the identity is in every cell because $S_1$ is inverse-closed.

Let $1 \leq i \leq \ell$, and let $t \in C_i$.
Since $S_1$ is inverse-closed, we can write $t = s_is_j^{-1}$ for some $1 \leq j \leq \ell$. Then
\[
t^{-1} C_i = \lsb t^{-1} s_i s_1, t^{-1} s_i s_2, \ldots, t^{-1} s_i s_\ell \rsb = \lsb s_js_1, s_js_2, \ldots, s_js_\ell \rsb = C_j \in \pi,
\]
so $\pi$ satisfies the \nameref{abc}.

Since $\Cay(G, S_1)$ is bipartite, by Proposition~\ref{cayleyBipartite}, we have a homomorphism $\varphi: G \rightarrow \mathbb{Z}_2$ such that $S_1 \subseteq \varphi^{-1}(1)$. We let $G_0= \ker(\varphi)$.
It follows from $S_1 \subseteq \varphi^{-1}(1)$ that $S_2 \subseteq \ker(\varphi) = G_0$. 

We now show that $\BCay(G_0, \pi) = \lp \bipartB \cup \bipartC, E \rp$ is isomorphic to $\Cay(G, S_1)$. We define $\psi: \bipartB \cup \bipartC \to G$ by
\begin{equation*}
\begin{cases}
h \mapsto h & h \in \bipartB.\\
hC_i \mapsto hs_i & hC_i \in \bipartC.
\end{cases}
\end{equation*} 
We first show that this map is well defined. This needs to be verified for vertices $x \in \bipartC$. Suppose that $h_1C_i = h_2C_j$. Then $h_1^{-1}h_2 \in C_i$ and $h_2^{-1}h_1 \in C_j$, but this is only possible if $h_1^{-1}h_2= s_is_j^{-1}$. Then $h_1s_i = h_2 s_j$, as desired.

Next, we show it is bijective. Suppose that $\psi(gC_i)=\psi(hC_j)$. Then, $gs_i=hs_j$ and thus $g^{-1}h =s_is_j^{-1} \in C_i$ and $h^{-1}g =s_js_i^{-1} \in C_j$. But this implies that $gC_i=hC_j$. This shows that $\psi$ is injective, and since both sets have the same cardinality, we have a bijection.

It is clear that the image of any edge in $\BCay(G_0, \pi)$ is an edge in $\Cay(G, S_1)$, and that any edge in $\Cay(G,S_1)$ is the image of some edge in  $\BCay(G_0, \pi)$, and hence $\psi$ is a graph isomorphism.\qed

\begin{rem}
We can also prove Theorem~\ref{cayleyBiCayley} using Theorem~\ref{hyperSabidussi}. Let $\sigma$ be as defined in Proposition~\ref{bipartite_auts}. Consider the hypergraph with vertex set $G_0 = \ker \lp \sigma \rp$ and hyperedge set $\lsb N \lp x \rp \subseteq G_0 : x \in \sigma^{-1} \lp 1 \rp \rsb$ where $N \lp x \rp$ is the neighbourhood of $x$. We can verify that this hypergraph has $\Cay(G,S)$ as an incidence graph. Since $G_0$ acts regularly on the vertex set of the hypergraph, Theorem~\ref{hyperSabidussi} implies that its incidence graph is a Cayley incidence graph.
\end{rem}

We now consider when a Cayley incidence graph is a Cayley graph. It clearly needs to be regular, but it is not the case that a regular Cayley incidence graph is always a Cayley graph, as evidenced by Example~\ref{exNontransitive}.
Below we give some sufficient conditions for when this is the case.

\begin{thm}\label{thmCayley} 
Let $\BCay(G, \pi)$ be a connected $\bipartC$-regular Cayley incidence graph with $\pi=\{C_1, C_2, \dots, C_k\}$,  $C_1= \{e, g_2, \ldots, g_k\}$, $g_1 = e$ and $C_i = g_i^{-1}C_1$. Let $C_i^* = g_iC_1^{-1}$ for $1\leq i\leq k$. If there exists a group automorphism $\varphi: G\rightarrow G$ of order $2$ such that $\varphi(C_1) = C_j^{*}$, then $\BCay(G, \pi)$ is a Cayley graph for the group $G' = G\rtimes \langle \varphi \rangle$. 
\end{thm}
\proof
The group automorphism $\varphi$ defines a graph isomorphism by Theorem \ref{Normalizer-duals}. It is easy to see that $G\rtimes \langle \varphi \rangle$ acts regularly on the vertices of $\BCay(G, \pi)$. Then by Sabidussi's theorem (Theorem~\ref{sabidussi}), we see that $\BCay(G, \pi)$ is a Cayley graph.\qed

For any abelian group $G$, the \textit{generalized dihedral group} $\dih(G)$ is the semidirect product $G \rtimes \Z_2$ with $ag = g^{-1}a$ for the non-identity element $a \in \Z_2$ and $g \in G$.

\begin{cor}\label{cor-gendih}
    Let $X = \BCay \lp G, \pi \rp$ be a connected Cayley incidence graph for an abelian group $G$. If $X$ is $\bipartC$-regular, then $X$ is a Cayley graph of $\dih(G)$.
\end{cor}
\proof There is a group automorphism on $G$ that maps $g$ to $g^{-1}$, so the result follows by Theorem~\ref{thmCayley}.\qed

\begin{cor}\label{thmZnCayley}A regular Cayley incidence graph of the cyclic group $\mathbb{Z}_n$ is a Cayley graph for the dihedral group of order $2n$.
\end{cor}
\proof By Lemma~\ref{lemAbellgeqk}, a regular $\BCay(G, \pi)$ with abelian $G$ is either $\bipartC$-transitive, or $\pi$ is a collection of $k$ subgroups of size $k$.
The latter is impossible for $G = \Z_n$, so Corollary~\ref{cor-gendih} applies.\qed

From Theorem~\ref{cayleyBiCayley} and Theorem~\ref{thmCayley}, we see that Cayley incidence graphs generalize bipartite Cayley graphs of girth at least six. An older generalization of Cayley graphs comes from bi-Cayley graphs.

\begin{defn}
    A \textit{bi-Cayley graph}, sometimes called a \textit{semi-Cayley graph}, is a graph which admits a semiregular group action such that there are exactly two vertex orbits of equal size.
\end{defn}

Bi-Cayley graphs are a natural way of loosening Sabidussi's thoerem (Theorem~\ref{sabidussi}) to include large classes of highly symmetric graphs. There has been much work on bi-Cayley graphs, from early results characterizing when bi-Cayley graphs were strongly regular~\cite{de1992strongly, malnivc2007strongly, maruvsivc1987strongly} to more recent work on the spectrum~\cite{arezoomand2015classification, gao2010spectrum}, particular classes of bi-Cayley graphs~\cite{conder2020edge, xie2024isomorphisms} and their properties~\cite{duan2023hamiltonian, wang2023perfect}. The work on bi-Cayley graphs is much more extensive than the examples listed here, but in the remainder of this section we are primarily concerned with comparing bi-Cayley graphs and their automorphisms to Cayley incidence graphs.

A bi-Cayley graph $\textup{BiCay}(G, R, L, S)$ is the graph with vertex set $G\times\{0,1\}$, and edge set
\begin{align*}
    & \{(g,0)\sim (gr,0) ~: ~ g\in G, r\in R\}\\
    &\cup \{(g,1)\sim (gl,1) ~: ~ g\in G, l\in L\}\\
    &\cup \{(g,0)\sim (gs,1) ~: ~ g\in G, s\in S\}.
\end{align*}
Thus any bipartite bi-Cayley graph such that the two orbits are the two biparts has the form $\textup{BiCay} \lp G, \emptyset, \emptyset, S \rp$.

It is easy to see that a $\bipartC$-regular Cayley incidence graph is a bi-Cayley graph. The converse also holds when the bi-Cayley graph is bipartite with girth at least six.

\begin{thm}\label{BicCayG00S}
Let $X= \textup{BiCay}(G, \emptyset, \emptyset, S)$, where $S=\{s_1, \dots, s_\ell\}$. If $X$ has girth at least six, then $X$ is isomorphic to $\BCay(G,\pi)$, where $\pi= \{C_1, \dots, C_\ell\}$, with
\begin{align*}
C_1 &=\{e, s_1^{-1}s_2,\dots, s_1^{-1}s_\ell\},\\
C_2 &=\{s_2^{-1}s_1, e,\dots, s_2^{-1}s_\ell \},\\
 & ~ ~ \vdots\\
C_\ell &= \{s_\ell^{-1} s_1, s_\ell^{-1} s_2 ,\dots, e\}.
\end{align*}

\end{thm}
\proof
It is easy to verify that $\pi$, satisfies the \nameref{abc}. We define a graph isomorphism from BiCay$\lp G, \emptyset, \emptyset, S \rp$ to $\BCay \lp G, \pi \rp$ by
\begin{equation*}
    \begin{cases}
(g, 0) \mapsto g,\\
(h, 1) \mapsto (hs_1) C_1.
\end{cases}
\end{equation*}
We now verify that this is an isomorphism. We see that $\{(g, 0), (h,1)\}$ is an edge in $\textup{BiCay}(G, \emptyset, \emptyset, S)$ if and only $h=gs$, where $s\in S$. On the other hand, we see that if $g$ and $(hs_1) C_1$ share an edge in $\BCay(G, \pi)$ if and only if $g=(hs_1)(s_1^{-1} s_j)$, for some $j$, i.e. if $g=hs_j$.\qed

It can be shown that Theorem~\ref{Normalizer-duals} is essentially equivalent to Theorem 1.1 from Zhou and Feng~\cite{Zhou2016} applied to bipartite bi-Cayley graphs, though we omit the details.

\section{Open Problems}

There is a well-developed theory of Cayley graphs, and it would be nice to develop a parallel theory for Cayley incidence graphs. We conclude by highlighting some open problems inspired by the work in the paper and the tables in the appendices.

Every group admits a Cayley graph, but the \nameref{abc} is more restrictive, and there are many groups that do not admit a non-trivial collection $\pi$ satisfying the conditions of Definition~\ref{bcayDef}, as summarized in Table~\ref{apAll}.

\begin{quest}
    When is a Cayley graph the underlying Cayley graph of a non-trivial Cayley incidence graph?
\end{quest}

Let $X = \Cay \lp G, S \rp$. For $X$ to be the underlying Cayley graph of a nontrivial Cayley incidence graph, we have to be able to partition the connection set into $\ell$ cells of size $k-1$, for $\ell \geq 2$ and $k \geq 3.$ Further, $X$ must contain  cliques of size $k$, and by Proposition~\ref{bCSpectrum}, the minimal eigenvalue of $X$ must be at least $-\ell.$ This imposes some necessary conditions, but as is clear from Table~\ref{apAll}, the structure of the group plays a vital role in determining how many $\pi$ satisfy the \nameref{abc}.

If a full characterization is out of reach, it would be helpful to be to build up Cayley incidence graphs as we did in Lemma~\ref{productCIN} and Lemma~\ref{lem:intersection}.

\begin{quest}
   What operations on groups and cells preserve the \nameref{abc}? 
\end{quest}

In the small examples listed in Table~\ref{apNonAbelian}, every regular Cayley incidence graph is isomorphic to a Cayley graph, but Example~\ref{exNontransitive} shows this need not be the case. Thus we might ask to build on Table~\ref{apNonAbelian} and the work in Section~\ref{secBiBi} to answer the following question.

\begin{quest}
    When is a regular Cayley incidence graph not a Cayley graph?
\end{quest}

It would also be interesting to further develop the connections between Cayley incidence graphs and the work in related areas, including hypergraphs, design theory, and bi-Cayley graphs.

Cayley hypergraphs have been introduced a number of times~\cite{buratti1994cayley, jajcayova2024generalizations, lee2013cayley, shee1990group}, and we can ask how the new class of Cayley incidence graphs can be used in relation to the problems considered by those authors. Of particular note, we highlight the biregular cages introduced in Section~\ref{secCages}.

\begin{prob}
    Use Cayley incidence graphs to find biregular cages.
\end{prob}

Section~\ref{secDifference} discussed the use of difference sets and partial difference families to construct projective and affine spaces as Cayley incidence graphs. 

\begin{quest}
    What other structures from design theory and finite geometry can be constructed using Cayley incidence graphs?
\end{quest}

A key property of Cayley graphs is that they give many examples of vertex-transitive graphs. This suggests that the $\bipartC$-transitive graphs introduced in Section~\ref{secIncMin} would provide an especially good extension of Cayley graphs to a bipartite context.

\begin{prob}
   Further develop the theory of $\bipartC$-transitive Cayley incidence graphs. 
\end{prob}

Finally, there is a beautiful theory of graphical regular representations of Cayley graphs, and it would be interesting to try and extend it to Cayley incidence graphs.

\begin{prob}
    When is the automorphism group of a Cayley incidence graph $\BCay \lp G, \pi \rp$ isomorphic to $G$?
\end{prob}

We could also reframe this problem in terms of hypergraphs.

\begin{prob}
    When is the automorphism group of the hypergraph with incidence graph $\BCay \lp G, \pi \rp$ isomorphic to $G$?
\end{prob}

\section*{Acknowledgements}
We would like to thank Jakob Hultgren, Maryam Sharifzadeh, Klara Stokes for organizing the DiGResS workshop at Ume\aa{} where this work was initiated. The first author is supported by CNPq, National Council for Scientific and Technological Development (Brazil). The second, third, and fourth authors are all supported by Kempe foundation JCSMK23-0058, JCSMK22-0160, and JCSMK21-0074 respectively. The fifth author is supported by the Wallenberg AI, Autonomous Systems and Software Program
(WASP) funded by the Knut and Alice Wallenberg Foundation.


\newpage
\appendix
\section*{Appendix}

\begin{table}[h!]
\begin{center}
\begin{tabular}{|c|c|c|c|}
    \hline
    Order & Group & Non-trivial $\pi$ \\ 
    \hline
    7 & $\Z_7$ & 1\\
    \hline
    8 & $\Z_8$ & 1\\
    \hline
    8 & $\Z_4\times \Z_2$  & 0\\
    \hline
    8 & $\Z_2^3$ & 0 \\  
    \hline
    8 & $D_4$ & 0\\
    \hline
    8 & $Q_8$ & 1\\
    \hline
    9 &$\Z_9$ & 1\\
    \hline
    9 & $\Z_3^2$ & 3\\
    \hline
    10 &$\Z_{10}$ & 1 \\
    \hline
    10 & $D_5$ & 0\\
    \hline
    11 &$\Z_{11}$ & 1 \\
    \hline
    12 &$\Z_{12}$ & 4\\
    \hline
    12 & $\Z_6\times\Z_2$ & 2\\
    \hline
    12 & $D_6$ & 0\\
    \hline
    12 & $\dic_{3}$ & 3\\
    \hline
    12 & $A_4$ & 3\\
    \hline
    13 &$\Z_{13}$ & 4\\
    \hline
    14&$\Z_{14}$ & 3\\
    \hline
    14 & $D_{7}$ & 0\\
    \hline
    15&$\Z_{15}$ & 16\\
    \hline
    16&$\Z_{16}$ & 8\\
    \hline
    16 & $\Z_4^2$ & 6\\
    \hline
    16 & $\Z_8\times\Z_2$ & 1\\
    \hline
    16 & $\Z_4\times\Z_2^2$ & 1\\
    \hline
    16 & $\Z_2^4$ & 4\\
    \hline
    16 & $D_8$ & 0\\
    \hline
    16 & $\dic_{4}$ & 5\\
    \hline
    16 & $\Z_4\rtimes \Z_4$ & 3\\
    \hline
    16 & $(\Z_2\times\Z_2)\rtimes \Z_4$ & 4\\
    \hline
    16 & $\Z_8\rtimes \Z_2$ & 2\\
    \hline
    16 & $QD_{8}$ & 2\\
    \hline
    16 & $D_4\times\Z_2$ & 2\\
    \hline
    16 & $Q_8\times\Z_2$ & 2\\
    \hline
    16 & $(\Z_4\times\Z_2)\rtimes\Z_2$ & 2\\
    \hline
\end{tabular}
\end{center}
\caption{List of all groups of order at most 16 (and at least seven) and the number of non-isomorphic, non-trivial Cayley incidence graphs for each. Here, $D_n$ is a dihedral group,  $\dic_n$ a dicyclic group, $Q_8$ the quaternion group, $A_4$ the alternating group, $QD_8$ the quasi-dihedral group, $Q_8 \times \Z_2$ the Hamiltonian group, and $(\Z_4\times\Z_2)\rtimes\Z_2$ the Pauli group.}
\label{apAll}
\end{table}

\newpage
\renewcommand{\arraystretch}{1.1}
\begin{table}[h!]
\begin{center}
\begin{tabular}{|c|c| c|c|c|c|c|c|}
\hline
 Group & $\ell$ & $k$ & Cayley & Order of $\Aut$& Orbits of $\Aut$ \\ 
 \hline
 $\Z_8$ &3&3& Yes &96&1 \\  
 \hline
 $\Z_9$ &3&3& Yes &18&1 \\ 
 \hline
 $\Z_{10}$ &3&3& Yes &20&1 \\  
 \hline
 $\Z_{11}$ &3&3& Yes &22&1 \\  
 \hline
 $\Z_{12}$ &3&3& Yes &24&1 \\
  \hline
 $\Z_{12}$ &3&3& Yes &24&1 \\
 \hline
 $\Z_{12}$ &3&3& Yes &48&1 \\
  \hline
  $\Z_{12}$ &4&3& No &12&3 \\
  \hline
  $\Z_{13}$ &3&3& Yes &26&1 \\
  \hline
 $\Z_{13}$ &3&3& Yes &78&1 \\
  \hline
  $\Z_{13}$ &6&3& No &39& 3 \\
  \hline
 $\Z_{13}$ &4&4& Yes &11232& 1 \\
  \hline
 $\Z_{14}$ &3&3& Yes &28&1  \\
  \hline
 $\Z_{14}$ &3&3& Yes &28&1 \\
  \hline
 $\Z_{14}$ &4&4& Yes &672&1 \\
  \hline
 $\Z_{15}$ &3&3& Yes &30&1 \\
  \hline
   $\Z_{15}$ &3&3& Yes &30&1 \\
  \hline
 $\Z_{15}$ &3&3& Yes &60&1 \\
  \hline
   $\Z_{15}$ &3&3& Yes &60&1 \\
  \hline
   $\Z_{15}$ &4&3& No &15&3 \\
  \hline
    $\Z_{15}$ &4&3& No &30&3 \\
  \hline
  $\Z_{15}$&6&3& No &15& 3 \\
  \hline
  $\Z_{15}$&6&3& No &15& 3 \\
  \hline
  $\Z_{15}$&6&3& No &15& 3 \\
  \hline
  $\Z_{15}$&6&3& No &360& 2 \\
  \hline
  $\Z_{15}$&6&3& No &60&2 \\
  \hline
  $\Z_{15}$&4&4& Yes &30& 1 \\
  \hline
    $\Z_{15}$&4&4& Yes &60& 1 \\
  \hline
    $\Z_{15}$&4&4& Yes &720& 1 \\
  \hline
    $\Z_{15}$&7&3& No &60& 3 \\
  \hline
    $\Z_{15}$&7&3& No &20160 = $8!/2$& 2 \\
  \hline
\end{tabular}
\end{center}
\caption{All non-trivial Cayley incidence graphs (up to graph isomorphism) for cyclic groups of order at most 15.}
\label{apCyclic}
\end{table}

\newpage


\begin{table}[h!]
\begin{center}
\begin{tabular}{|c|c|c|c|c|c|}
\hline
 Group& $\ell$ & $k$ & Cayley & Order of $\Aut$& Orbits of $\Aut$ \\ 
 \hline
 $\Z_3^2$ &2&3& No &72&2 \\
 \hline 
  $\Z_3^2$ &3&3& Yes &216&1 \\
 \hline 
 $\Z_3^2$ &4&3& No &432& 2 \\
 \hline 
  $\Z_6\times\Z_2$ &3&3& Yes &144& 1 \\
 \hline 
   $\Z_6\times\Z_2$ &4&3& No &576& 2 \\
  \hline 
   $\Z_4^2$ &3&3& Yes &192& 1 \\
 \hline 
$\Z_4^2$ &2&4& No &1152& 2 \\
 \hline 
 $\Z_4^2$ &3&4& No &192& 2 \\
 \hline 
  $\Z_4^2$ &6&3& No &96& 2 \\
 \hline 
   $\Z_4^2$ &4&4& Yes &2304& 1 \\
 \hline 
    $\Z_4^2$ &5&4& No &5760& 2 \\
 \hline 
    $\Z_8\times\Z_2$ &3&3& Yes &64& 1 \\
 \hline 
    $\Z_4\times\Z_2\times\Z_2$ &2&4& No &1152&2 \\
 \hline
    $\Z_2^4$ &2&4& No &1152&2 \\
 \hline
    $\Z_2^4$ &3&4& No &576&2 \\
 \hline
    $\Z_2^4$ &4&4& Yes &2304&1 \\
 \hline
     $\Z_2^4$ &5&4& No &5760&2 \\
 \hline
\end{tabular}
\end{center}
\caption{All non-trivial Cayley incidence graphs (up to graph isomorphism) for abelian, non-cyclic groups of order at most 16.}
\label{apAbelian}
\end{table}

\newpage


\begin{table}[h!]
\begin{center}
\begin{tabular}{|c|c|c|c|c|c|c|}
\hline
 Group & $\ell$ & $k$ & Girth & Cayley & Order of $\Aut$& Orbits of $\Aut$ \\ 
 \hline
 $Q_8$ & 3&3&6& Yes &96&1 \\  
  \hline
 $\dic_{12}$ &3&3&6&Yes&48&1\\  
 \hline
 $\dic_{12}$ &3&3&6&Yes&144&1\\  
 \hline
$\dic_{12}$ &4&3&6&No&576 & 2 \\  
 \hline
 $A_4$ & 3&3&6&Yes&144&1 \\  
 \hline
  $A_4$ &2&4&6&No&48&2  \\  
 \hline
 $A_4$ &4&3&6&No&576&2 \\  
 \hline
 $\dic_4$ &3&3&6&Yes&64&1 \\  
 \hline
 $\dic_4$ &6&3&6&No&32&3 \\  
 \hline
 $\dic_4$ &4&4&6&Yes&2304&1\\  
 \hline
 $\dic_4$ &4&4&6&Yes&64&1\\  
 \hline
 $\dic_4$ &5&4&6&No&5760&2\\  
 \hline
 $\Z_4\rtimes \Z_4$ &3&3&6&Yes&64&1\\  
 \hline
 $\Z_4\rtimes \Z_4$ &2&4&8&No&1152&2 \\  
 \hline
  $\Z_4\rtimes \Z_4$ &6&3&6&No&96&2 \\  
 \hline
 $\Z_2^2\rtimes \Z_4$ &2&4&8&No&1152&2\\  
 \hline
 $\Z_2^2\rtimes \Z_4$ &3&4&6&No&576&2 \\  
 \hline
 $\Z_2^2\rtimes \Z_4$ &4&4&6&Yes&2304&1 \\  
 \hline
 $\Z_2^2\rtimes \Z_4$ &5&4&6&No&5760&2 \\  
 \hline
 $\Z_8\times\Z_2$ &2&4&8&No&1152&2\\  
 \hline
  $\Z_8\times\Z_2$ &3&4&6&No&192&2 \\  
 \hline
 $QD_8$ &2&4&8&No&1152&2\\  
 \hline
 $QD_8$ &3&4&6&No&192&2\\  
 \hline
 $D_4\times \Z_2$ &2&4&8&No&1152&2\\  
 \hline
 $D_4\times \Z_2$ &3&4&6&No&192&2 \\  
 \hline
 $Q_8\times \Z_2$ &4&4&6&Yes&2304&1\\  
 \hline
 $Q_8\times \Z_2$ &5&4&6&No&5760&2\\  
 \hline
 $(\Z_4\times\Z_2)\rtimes\Z_2$ &2&4&8&No&1152&2\\  
 \hline
 $(\Z_4\times\Z_2)\rtimes\Z_2$ &3&4&6&No&192&2\\
 \hline
 \end{tabular}
 \end{center}
\caption{All non-trivial Cayley incidence graphs (up to graph isomorphism) for non-abelian groups of order at most 16 with $D_n$ a dihedral group,  $\dic_n$ a dicyclic group, $Q_8$ the quaternion group, $A_4$ the alternating group, $QD_8$ the quasi-dihedral group, $Q_8 \times \Z_2$ the Hamiltonian group, and $(\Z_4\times\Z_2)\rtimes\Z_2$ the Pauli group.}
\label{apNonAbelian}
\end{table}

\end{document}